%% file: main.tex
\documentclass[a4paper,11pt]{amsart}
\pdfoutput=1
\input{preamble.tex}

\usepackage{xstring}

\newcommand{\errorplot}[4][]{%
    \IfBeginWith{#2}{qsvt}{\def\stepmult{2}}{%
        \IfBeginWith{#2}{q_cheb}{\def\stepmult{2}}{%
        \IfBeginWith{#2}{chebopt}{\def\stepmult{2}}{\def\stepmult{1}}%
        }%
    }%
    \IfBeginWith{#2}{qsvt}{\def\stepadd{2}}{%
        \IfBeginWith{#2}{q_cheb}{\def\stepadd{2}}{%
        \IfBeginWith{#2}{chebopt}{\def\stepadd{2}}{\def\stepadd{1}}%
        }%
    }%
    \addplot+[
        #1,unbounded coords=discard,mark options={solid},
    ] table [
        col sep=semicolon,
        x expr={(\thisrow{steps} * \stepmult + \stepadd)},
        y={#2_50}
    ] {#3};
    \addplot+[
        #1,no marks,dotted,unbounded coords=discard,mark options={solid},forget plot
    ] table [
        col sep=semicolon,
        x expr={(\thisrow{steps} * \stepmult + \stepadd)},
        y={#2_95}
    ] {#3};
}

\newcommand{\complexityplot}[3][]{
\addplot+[
    #1,unbounded coords=discard,mark options={solid}
] table [
    col sep=semicolon,
    x={#2_complexity},
    y={#2_50},
] {#3};
}

\definecolor{color1}{RGB}{64, 83, 211}
\definecolor{color2}{RGB}{221, 179, 16}
\definecolor{color3}{RGB}{181, 29, 20}
\definecolor{color4}{RGB}{0, 190, 255}
\definecolor{color5}{RGB}{251, 73, 176}
\definecolor{color6}{RGB}{0, 178, 93}
\definecolor{superlightgray}{RGB}{240, 240, 240}

\newcommand{\N}{\mathbb{N}}

\newcommand{\R}{\mathbb{R}}
\newcommand{\C}{\mathbb{C}}
\newcommand{\E}{\mathbb{E}}

\newcommand{\calS}{\mathcal{S}}

\newcommand{\calO}{\mathcal{O}}

\DeclareMathOperator{\diag}{diag}
\DeclareMathOperator*{\argmin}{arg\,min}
\DeclareMathOperator*{\argmax}{arg\,max}

\newcommand{\Id}{\mathrm{Id}}

\newcommand{\caps}{\mathrm{CAP}}
\newcommand{\cups}{\mathrm{CUP}}
\newcommand{\qsvt}{\mathrm{QSVT}}
\newcommand{\opt}{\mathrm{opt}}
\newcommand{\cheb}{\mathrm{Cheb1}}
\newcommand{\qcheb}{\mathrm{Cheb2}}
\newcommand{\chebopt}{\mathrm{Cheb3}}
\newcommand{\err}{\mathrm{err}}
\newcommand{\samples}{N}
\newcommand{\Nquadrature}{N_\mathrm{q}}
\newcommand{\Nnoise}{N_\mathrm{noise}}

\newcommand{\tol}{\varepsilon}

\begin{document}

\title[Constrained Optimal Polynomials for Quantum Linear Solvers]{Constrained Optimal Polynomials for Quantum Linear System Solvers}
\author[Matthias Deiml]{Matthias Deiml${}^\dag$}
\author[Daniel Peterseim]{Daniel Peterseim${}^\ddagger$}
\address{${}^\dag$ Institute of Mathematics, University of Augsburg, Universit\"atsstr.~12a, 86159 Augsburg, Germany}
\address{${}^\ddagger$ Institute of Mathematics \& Centre for Advanced Analytics and Predictive Sciences (CAAPS), University of Augsburg, Universit\"atsstr.~12a, 86159 Augsburg, Germany}
\thanks{Funded by the Deutsche Forschungsgemeinschaft (DFG, German Research Foundation) -- 571768116. We also acknowledge the use of IBM Quantum Credits for this work. The
views expressed are those of the authors, and do not reflect the official
policy or position of IBM or the IBM Quantum team.}
\date{April 28, 2026}
\email{\{matthias.deiml,daniel.peterseim\}@uni-a.de}

\markboth{\textls{\textsc{\MakeLowercase{\shortauthors}}}}{\textls{\textsc{\MakeLowercase{\shorttitle}}}} % set running head
% 68Q12, 65F10, 81P68, 65N22 

\begin{abstract}
Quantum linear system solvers typically realize the inverse map as a polynomial transformation of the spectrum, so their practical cost hinges on implementing this transformation at a low polynomial degree. We introduce constrained optimal polynomials as a framework for this task, drawing on classical Krylov subspace theory. Within this framework, we develop two classes of solvers. Constrained Uniform Polynomial (CUP) solvers optimize the tradeoff between approximation accuracy and block encoding normalization under a uniform spectral model consistent with the available bounds. Constrained Adaptive Polynomial (CAP) solvers retain this structure but replace the uniform model with a probability measure reconstructed from spectral moments via a maximum entropy ansatz, where the moments are extracted from QSVT measurements. Numerical experiments under hardware and stochastic noise show that these methods achieve lower error than standard QSVT-based and Chebyshev-iteration-type solvers, particularly in noise-limited regimes. CUP offers robust performance under generic spectra, while CAP provides further improvement when the spectral structure can be exploited.

\item[\hskip\labelsep\scshape{}Keywords.] quantum linear system solver, quantum singular value transformation, polynomial approximation, Krylov subspace method, block encoding, maximum entropy method, noisy quantum computation
\end{abstract}

\maketitle

\section{Introduction}
Solving linear systems is a fundamental problem in scientific computing and has attracted significant attention in the quantum setting since the seminal work~\cite{HHL09}. Various quantum algorithms have since been proposed, including improvements based on quantum signal processing and block-encoding techniques~\cite{CKS17,CGJ19,GSLW19}. An overview is given in~\cite{MPS+25}. These algorithms exploit the ability of quantum computers to represent and transform vectors in a way that can yield a logarithmic dependence on the dimension $d$. When the matrix admits an efficient quantum implementation, the best of these algorithms, based on Variable Time Amplitude Amplification (VTAA)~\cite{Amb10,LS26} or adiabatic techniques and eigenfiltering~\cite{LT20,AL22,Dal24}, achieve an asymptotic runtime bound that scales~as
\[
\calO(\tol^{-1} \kappa \log(d\tol^{-1})),
\]
where $\tol$ is the error tolerance and $\kappa$ is the condition number of the linear system. This is to be compared with the classical iterative cost of $\calO(d\kappa \log \tol^{-1})$ when matrix--vector multiplication is linear in $d$. The quantum bound improves the dependence on $d$ from linear to logarithmic, but degrades the dependence on $\tol^{-1}$ from logarithmic to linear.
 In most leading quantum linear solvers, the inverse map is realized through a polynomial transformation of the spectrum implemented via quantum singular value transformation (QSVT)~\cite{GSLW19} (see also \cref{sec:qsvt}). Practical performance is therefore dictated by the accuracy of this transformation at a given polynomial degree. Two closely related issues are particularly important.
\begin{enumerate}
\item The accuracy of polynomial approximations at a pre-asymptotic, low polynomial degree is a central practical issue since the degree directly controls circuit depth. Classical numerical linear algebra offers a broad range of techniques for optimal polynomial approximation, many of which are not reflected in current quantum methods. In particular, insights from nearly a century of research on Krylov subspace methods remain largely unexplored in the design of quantum linear system solvers. Related concerns already appear in the classical literature; for example, \cite{FS50} discusses limitations of approximations of the type used in~\cite[Theorem~41]{GSLW19}; see also \cite[Section~5.5.3]{LS12} for historical context. Only recently were some of these concerns addressed in the quantum context in \cite{GKS24,SNW+25} through the definition of Chebyshev-iteration-type solvers (\cref{sec:chebsym}).
\item Because hardware noise and measurement errors dominate on current and foreseeable devices, theoretical runtime bounds are not a reliable guide to practical performance. Even recent static~\cite{GKS24,SNW+25} and adaptive approaches~\cite{TWYH24,XZZ24,LWX25} do not systematically analyze performance in the presence of hardware noise or measurement errors. Sup-norm reduction is at best a post-processing step applied to a polynomial constructed for approximation accuracy, rather than part of a joint optimization.
\end{enumerate}
Motivated by these considerations, we develop two classes of quantum linear system solvers inspired by classical Krylov subspace techniques (see the textbooks~\cite{Saa03,LS12} for a comprehensive overview) and study their performance in the presence of hardware and stochastic noise. The first class comprises Constrained Uniform Polynomial (CUP) solvers (\cref{sec:cup}), which optimize the tradeoff between approximation accuracy and block encoding normalization under a uniform spectral model consistent with the available spectral bounds. The second class includes Constrained Adaptive Polynomial (CAP) solvers (\cref{sec:cap}), which retain the same constrained optimization but replace the uniform model with a problem-specific probability measure reconstructed from spectral moments via a~maximum entropy ansatz~\cite{MP84}, with moments extracted from QSVT measurements. CUP therefore captures the best uniform polynomial for the given spectral interval, while CAP adapts to the actual spectral distribution of the problem at hand. These solver polynomials can be combined with structural transformations of the block encoding domain (\cref{sec:transforms}) that further reduce the effective normalization interval when a block encoding of a matrix square root is available. 

Across benchmark computations (\cref{sec:num}), these methods outperform standard QSVT-based inversion as well as Chebyshev-iteration-type solvers under hardware and stochastic noise. The advantage is most pronounced in noise-limited regimes but persists in noiseless settings as well. 

The value of reducing the polynomial degree required to achieve a prescribed accuracy extends beyond near-term devices. In the fault-tolerant era, circuit depth will remain a central practical resource since longer coherent computations generally entail greater error correction overhead. Even in idealized settings, low-degree solver polynomials remain valuable when the linear system computation appears as a subroutine within a larger algorithm~\cite{DP24,BDP26}. In such contexts, realizing nonlinear transformations of the solution state may require multiple independent preparations, and the No-Cloning Theorem~\cite{WZ82} implies that these preparations must be carried out separately rather than obtained by duplication. The design of low-degree, efficiently implementable solver polynomials is, therefore, a central concern in both near-term and fault-tolerant regimes.

\section{Model problem}
Consider the linear system of equations
\begin{equation*} \label{eq:system}
Ax = b
\end{equation*}
with coefficient matrix~$A \in \R^{d \times d}$, right-hand side vector~$b \in \R^d$, and dimension~$d \in \N$. 
We focus on the case where $A$ is \emph{symmetric} and \emph{positive definite}. 
While several of the constructions considered in this paper can also be extended, with suitable modifications, to non-symmetric and indefinite matrices, the positive definite setting is already of central interest in its own right, for example, for discretizations of elliptic partial differential equations. It also provides a particularly natural setting for introducing the main ideas underlying constrained polynomial design.

We assume that the input data are provided via so-called \emph{block encodings}. For the purposes of this paper, this means that a matrix or vector~$B$ is realized by a quantum circuit with two associated parameters: the \emph{gate count}~$T_B \in \N$, i.e., the number of elementary gates required by the circuit, and the \emph{normalization}~$\gamma_B \ge \|B\|$, a scaling factor chosen so that $B/\gamma_B$ can be embedded as a subblock of a unitary. Lower values of $T_B$ and $\gamma_B$ correspond to a more efficient block encoding, with $\gamma_B$ entering the runtime of subsequent quantum subroutines polynomially. 
Here, $\|B\|$ denotes the spectral norm, i.e., the matrix norm induced by the Euclidean vector norm. For column vectors, this coincides with the Euclidean norm, and for symmetric positive definite matrices, it equals the largest eigenvalue. For the formal definition and further discussion of block encodings, we refer to~\cite{GSLW19,DP25}.

As an intermediate step, most quantum solvers produce a block encoding of a vector that approximates the unknown solution $x \in \R^d$. This can, in principle, be done efficiently since quantum computers can apply certain matrix--vector products using only a logarithmic number of gates $T_A \in \calO(\log d)$, relative to the dimension~$d$. To retain this potential advantage, it is essential to avoid recovering the full solution vector $x \in \R^d$, as this would generally require reading all $d$ entries and hence $\calO(d)$ gates. 

Instead, one typically seeks to extract only specific information from the encoded solution. In this work, we focus on quantities of interest of the form $x^\top M x$ for a given matrix~$M \in \R^{d\times d}$. Such quantities can be estimated, for example, by combining the Hadamard test~\cite{AJL06} with Monte Carlo sampling (empirical averaging of independent measurement outcomes). The following lemma states the resulting complexity bound. 
\begin{lemma}[Estimation procedure] \label{lem:estimation}
    Given block encodings of a vector~$x$ and a matrix~$M$, as well as an error tolerance $\tol > 0$ and failure probability~$0 < \delta < 1$, there is a procedure that uses $\samples(T_x + T_M)$ gates, where
    \[\samples \in \calO(\gamma_x^4\gamma_M^2 \tol^{-2} \log \delta^{-1}), \]
    and, with a probability of at least $1 - \delta$, returns an approximation of $x^\top Mx$ with an absolute error of at most $\tol$.
\end{lemma}
The proof follows from combining the Hadamard test with Hoeffding's inequality. In principle, there exist improved estimation algorithms, such as~\cite{SUR+20}, whose complexity depends only linearly on $\tol^{-1}$. However, these improvements are unlikely to be relevant on current hardware. More precisely, to obtain the improved linear dependence, the $\tol^{-1}T_x$ gates must be executed coherently, while Monte Carlo sampling uses $\tol^{-2}$ independent runs of only $T_x$ gates each. Even this bound is difficult to realize in practice. The dominant obstacle is hardware noise: the probability of a~fault-free execution decays with the circuit depth  $T_x + T_M$, so the achievable accuracy is ultimately limited by the circuit used to prepare the state being measured. Because theoretical complexity bounds do not adequately reflect this effect, we use the following error measure as the main benchmark for quantum linear system solvers.
\begin{definition}[Target quantities for quantum linear system solvers] \label{def:complexity}
For the purposes of this paper, a \emph{quantum linear system solver} is an algorithm~$\calS$ which, given block encodings of matrices $A$ and $M$ and a vector $b$, produces an estimate
    \[
    \calS(A, b, M) \approx b^\top A^{-1}MA^{-1}b=x^\top Mx.
    \]
    The \emph{complexity} of the solver is the number of gates it executes, excluding the cost of accessing $M$.
    The \emph{error} of the solver for fixed $A$ and $b$ is defined as
    \[\err(A, b) \coloneq \E\bigl[\bigl(\sqrt{\max\{0, \calS(A, b, \Pi_x)\}} - \|x\|\bigr)^2 + |\calS(A, b, \Id - \Pi_x)|\bigr] / \|x\|^2, \]
    where $\Pi_x \coloneq xx^\top/\|x\|^2$ is the projection onto $\operatorname{span}\{x\}$, and the expectation is taken over the sampling and hardware noise of $\calS$.
\end{definition}
Note that evaluating $\err(A, b)$  requires knowledge of the true solution $x$. It is used only as a~benchmark in our numerical experiments and not by the solver itself. The motivation for this error measure is as follows. If the solver is interpreted as producing an approximate solution vector~$\tilde x$, then the two target quantities correspond to the squared norms of the projections of~$\tilde x$ onto $\operatorname{span}\{x\}$ and its orthogonal complement:
\[
\calS(A, b, \Pi_x) \approx \|\Pi_x \tilde x\|^2, \qquad \calS(A, b, \Id - \Pi_x) \approx \|(\Id - \Pi_x) \tilde x\|^2.
\]
If these quantities were evaluated exactly, then $\err(A, b) = \|x - \tilde x\|^2 / \|x\|^2$. Moreover, if the approximation is exact, i.e.\ $x = \tilde x$, but the measurements are inexact, then the resulting measurement error is also reflected in $\err(A, b)$. The regularizations $\max\{0, \cdot\}$ and $|\cdot|$ ensure that $\err(A,b)$ is well-defined and nonnegative even when the noisy estimates of $\|\Pi_x \tilde x\|^2$ and $\|(\Id - \Pi_x)\tilde x\|^2$ take slightly negative values. Violations of nonnegativity contribute to the error rather than canceling~it.

\section{Quantum singular value transformation}\label{sec:qsvt}
When the dimension $d$ is large and $A$ is accessed only through matrix--vector products, both classical and quantum linear solvers are naturally restricted to approximations built from such products. In such settings, polynomial approximations based on repeated applications of~$A$ to~$b$ become a natural class of candidate solvers. Accordingly, one considers approximate solution vectors of the form
\begin{equation} \label{eq:krylov-approximation}
x \approx P(A)b = c_0 b + c_1 Ab + \dots + c_m A^m b,
\end{equation}
where $P$ is a polynomial of degree $\deg P = m$ with coefficients $c_0, \dots, c_m \in \R$.

Spectrally, the matrix polynomial~$P(A)$ can be written as
\[
P(A) = V \diag(P(\lambda_0), \dots, P(\lambda_{d-1})) V^\top,
\qquad\text{where}\qquad
A = V \diag(\lambda_0, \dots, \lambda_{d-1}) V^\top.
\]
Likewise,
\[
A^{-1} = V \diag(1/\lambda_0, \dots, 1/\lambda_{d-1}) V^\top.
\]
This shows that polynomial approximations to~$1/y$ on the spectrum of~$A$ yield polynomial approximations to the inverse~$A^{-1}$ in~\eqref{eq:krylov-approximation}. 

On a classical computer, the specific way in which matrix--vector multiplications are used to evaluate a polynomial is usually of secondary importance. Although certain factorizations are often preferred over the fully expanded form~\eqref{eq:krylov-approximation} for reasons of numerical stability, these differences disappear in exact arithmetic. On a quantum computer, by contrast, different implementations of the same polynomial can lead to substantially different normalizations. Since normalization directly affects the complexity limits in \cref{lem:estimation,def:complexity}, we seek implementations whose normalization is as small as possible. The best normalization that one can hope for is
\[
\gamma_{P(A)} = \|P\|_{[-\gamma_A, \gamma_A]} \coloneq \max_{y \in [-\gamma_A, \gamma_A]} |P(y)|.
\]
Quantum Singular Value Transformation (QSVT)~\cite{GSLW19} realizes this optimum up to a factor of $2$.

\begin{lemma}[QSVT] \label{lem:qsvt}
Let $A \in \R^{d \times d}$ be a symmetric matrix and $b \in \R^d$. Consider a real polynomial~$P$ with even and odd parts
\[P(y) = \underbrace{c_0 + c_2 y^2 + \dots}_{P_\mathrm{even}(y)} + \underbrace{c_1 y + c_3 y^3 + \dots}_{P_\mathrm{odd}(y)}.\]
Given block encodings of $A$ and $b$, we can obtain a block encoding of $P(A)b$
with a gate count scaling as \mbox{$T_{P(A)b} \in \calO(T_b + T_A \deg P)$} and a normalization~$\gamma_{P(A)b}$ with
\[\gamma_b \|P\|_{[-\gamma_A, \gamma_A]}\leq\gamma_{P(A)b} = \gamma_b (\|P_\mathrm{even}\|_{[-\gamma_A, \gamma_A]} + \|P_\mathrm{odd}\|_{[-\gamma_A, \gamma_A]}) \le 2 \gamma_b \|P\|_{[-\gamma_A, \gamma_A]}.\]
\end{lemma}

For specifics see \cite[Corollary~18]{GSLW19}. An important application is the solver~\cite[Theorem~41]{GSLW19} based on the analysis of~\cite[Theorem~4]{CKS17}, which we will refer to simply as the \emph{QSVT solver}. It is obtained by applying \cref{lem:qsvt} to the polynomial
\[
P_{\qsvt}^{(2n-1)}(y) \coloneq 4 \sum_{j = 0}^{n-1} (-1)^j\left[\frac{\sum_{k = j+1}^{\tilde{n}} \left(\begin{smallmatrix}2\tilde{n} \\ \tilde{n} + k\end{smallmatrix}\right)
}{2^{2\tilde{n}}} \right] T_{2j+1}(y),
\]
where $T_j$ is the $j$-th \emph{Chebyshev polynomial of the first kind} given by
\begin{equation} \label{eq:chebyshev}
T_j(y) = \cos(j \arccos(y)) \text{ for }y \in [-1, 1].
\end{equation}
The polynomial~$P_\qsvt$ approximates $1/y$ well on the intervals $[-1, -1/\kappa]$ and $[1/\kappa, 1]$ if $\tilde{n}$ is chosen as $\tilde{n} \coloneq \lceil \kappa^2 \log(\kappa / \tol) \rceil$, where $\tol$ is the measurement accuracy. This can be seen by setting $n = \tilde{n}$, in which case
\begin{equation}\label{e:Pqsvt}
P_{\qsvt}^{(2n-1)}(y) = \frac{1 - (1 - y^2)^n}{y},
\end{equation}
and noting that $1 - (1 - y^2)^n$ approaches $1$ for $y \neq 0$ and $n \to \infty$.
The parameter~$n \leq \tilde n$, which may be interpreted as the number of \emph{steps}, controls both the accuracy of the approximation and the degree of the polynomial. In the theoretical bounds, the required number of steps scales like~$\sqrt{\tilde n}$ and is therefore linear in the condition number~$\kappa$. Thus, at the level of polynomial degree, the standard QSVT construction already exhibits the desired linear dependence on~$\kappa$.

This does not, however, imply an overall solver complexity proportional to~$\kappa$. The cost of a~quantum linear solver depends not only on the polynomial degree but also on the normalization of the block encoding of the approximate solution. For the standard QSVT-based solver of~\cite[Theorem~41]{GSLW19}, this normalization may exceed the norm of the solution itself by a factor of~$\kappa$. Combined with the linear dependence of the polynomial degree on~$\kappa$, this yields an overall complexity proportional to~$\kappa^2$.

Optimal linear dependence on~$\kappa$ can be recovered by approaches based on variable-time techniques~\cite{Amb10,LS26}, as well as by adiabatic or eigenfiltering methods~\cite{LT20,AL22,Dal24}. Our focus in this article is instead on polynomial solver design within a common QSVT-based implementation framework. From a practical perspective, methods outside this framework may also incur a noticeable constant overhead compared with the standard QSVT solver, which can be disadvantageous in noise-limited regimes.

\section{Quantum Chebyshev iteration} \label{sec:chebsym}

\begin{figure}
    \centering
    \begin{tikzpicture}
    \begin{axis}[
        width=0.6\textwidth,
        height=0.28\textwidth,
        xlabel={$y$},
        domain=0:1,
        xmin=0, xmax=1.1, 
        ymin=0, ymax=10.9,
        samples=200,
        axis x line=center,
        axis y line=center,
        scale only axis,
        clip bounding box=upper bound,
        clip=true,
        font=\footnotesize,
        legend cell align=left,
        legend style={at={(1,0.95)},anchor=north east}
    ]
    
    \addplot[color=color6, thick] 
        {2.0052088 - 2.02071098 * (-1.0 + 2.0*x) + 1.82302404 * (-1.0 + 2.0*x)^2 - 1.42814673 * (-1.0 + 2.0*x)^3 + 4.03683881 * (-1.0 + 2.0*x)^4 - 7.24278805 * (-1.0 + 2.0*x)^5 - 7.39361727 * (-1.0 + 2.0*x)^6 + 16.66661552 * (-1.0 + 2.0*x)^7 + 16.90985939 * (-1.0 + 2.0*x)^8 - 22.7235355 * (-1.0 + 2.0*x)^9 };
    \addlegendentry{$m = d - 1$}
    \addplot[color=color5, thick] 
        {0.84251701 - 4.61424569 * (-1.0 + 2.0*x) + 6.08766338 * (-1.0 + 2.0*x)^2};
    \addlegendentry{$m = 2$}

    \addplot[color=black, dashed] 
        {1/x};
    \addlegendentry{$1/y$}

    \addplot+[
        scatter,
        color=black,
        mark=o,
        only marks,
        scatter src=explicit,
        scatter/use mapped color={fill=black,draw=black},
        scatter/@pre marker code/.append style=
                {/tikz/mark size={.6pt+\pgfplotspointmetatransformed/650}}
    ] table[x=x, y=y, meta=z, row sep=crcr] {
        x y z \\
        0.1 10. 0.02520761 \\
        0.29225254 3.42169826 0.61351211 \\
        0.32442912 3.08233736 0.20120367 \\
        0.36594263 2.73266881 0.23955599 \\
        0.59129006 1.69121734 0.03233618 \\
        0.64570662 1.54869096 0.12043292 \\
        0.69792337 1.43282207 0.13253725 \\
        0.12730387 7.85522063 0.64967664 \\
        0.73883432 1.35348343 0.18856219 \\
        0.98839558 1.01174066 0.18532726 \\
    };
    \addlegendentry{$\lambda_0, \dots, \lambda_{d-1}$}
    
    \draw node[color=color6, pin=240:{$\lambda_2$}] at (axis cs:0.29225254, 3.42169826) {};
    \end{axis}
\end{tikzpicture}
    \caption{Optimal polynomials with $d = 10$. The $y$ coordinate of the dots corresponds to the eigenvalues, their size to the weight $(V^\top b)_j$. The polynomial for $m=9$ interpolates at the eigenvalues. The polynomial for $m = 2$ approximates the extremal eigenvalues closely, as well as the eigenvalue $\lambda_2$ having large weight.}
    \label{fig:optimal}
\end{figure}
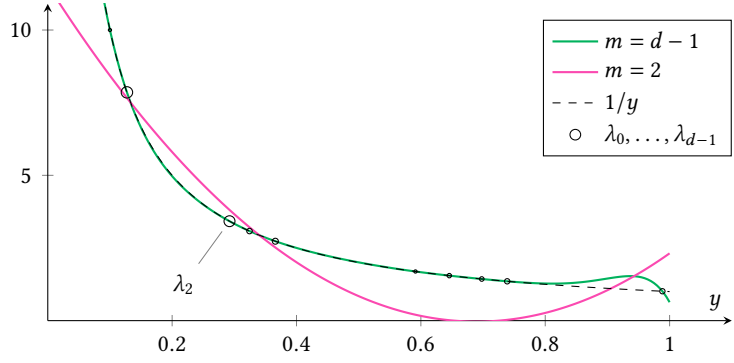

Although the QSVT solver is an efficient quantum linear system algorithm, its underlying polynomial $P_\qsvt$ is not optimal from the viewpoint of Krylov approximation.
Concretely, the best approximation $\smash{P_{\opt}}$ of the form~\eqref{eq:krylov-approximation} is characterized by
\begin{equation} \label{eq:argmin}
P_{\opt} =\argmin_{\deg P \le m}\|P(A)b - A^{-1}b\|^2 = \argmin_{\deg P \le m}\sum_{j=0}^{d-1} (P(\lambda_j) - \lambda_j^{-1})^2(V^\top b)_j^2.
\end{equation}
When $m = d - 1$, $\smash{P_{\opt}}$ interpolates $1/y$ at the eigenvalues of $A$, making the inverse exact; see also \cref{fig:optimal}. Of course, these eigenvalues are not known, but many solver polynomials effectively correspond to interpolations at suitable approximations of the eigenvalues.
Since the numerator $1 - (1 - y^2)^n$ in \eqref{e:Pqsvt} is strictly smaller than $1$ for $|y| < 1$, $P_\qsvt$ interpolates $1/y$ only at $y = \pm 1$; at least if $n = \tilde{n}$. Thus, there is no matrix such that this polynomial is optimal in the sense of the optimality condition~\eqref{eq:argmin}.

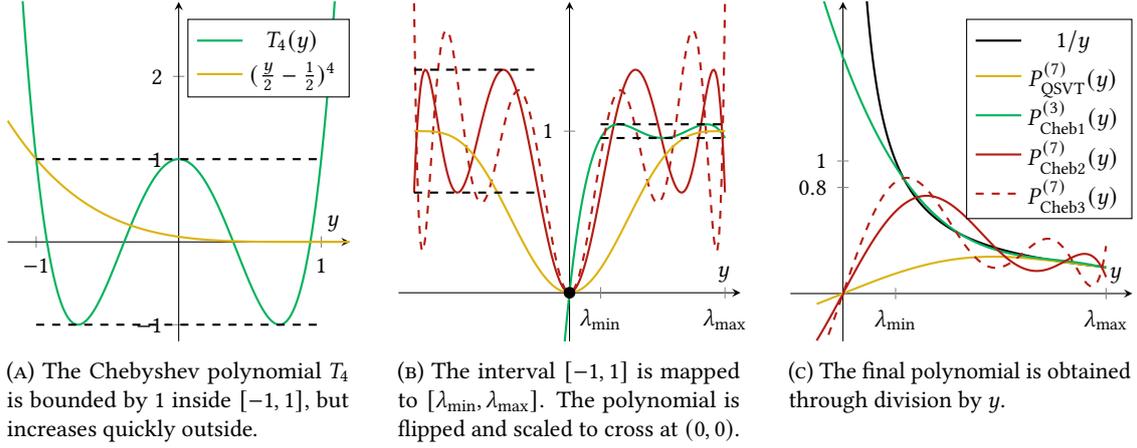
\begin{figure}
\begin{subfigure}[c]{0.3\textwidth}
    \begin{tikzpicture}
    \begin{axis}[
        width=\textwidth,
        height=\textwidth,
        xlabel={$y$},
        domain=-1.2:1.2,
        xmin=-1.2, xmax=1.2, 
        ymin=-1.2, ymax=2.9,
        samples=200,
        axis x line=center,
        axis y line=center,
        scale only axis,
        clip bounding box=upper bound,
        clip=true,
        font=\footnotesize,
        legend style={at={(1,0.95)},anchor=north east}
    ]
    
    \addplot[color=color6, thick] 
        {8 * x^4 - 8 * x^2 + 1};
    \addlegendentry{$T_4(y)$}
    \addplot[color=color2, thick] 
        {(x/2 - 1/2)^4};
    \addlegendentry{$(\tfrac{y}{2} - \tfrac12)^4$}
    \addplot[color=black, thick, dashed, domain=-1:1] {-1};
    \addplot[color=black, thick, dashed, domain=-1:1] {1};

    \end{axis}
\end{tikzpicture}
\subcaption{
The Chebyshev polynomial $T_4$ is bounded by $1$ inside $[-1, 1]$, but increases quickly outside.
}
\end{subfigure}
\hspace{0.03\textwidth}
\begin{subfigure}[c]{0.3\textwidth}
    \begin{tikzpicture}
    \begin{axis}[
        width=\textwidth,
        height=\textwidth,
        xlabel={$y$},
        domain=-5:5,
        xtick=\empty,
        ytick=\empty,
        extra x ticks={1,5},
        extra x tick labels={$\lambda_{\min}$, $\lambda_{\max}$},
        extra y ticks={1},
        xmin=-5.5, xmax=5.5, 
        ymin=-0.3, ymax=1.8,
        samples=200,
        axis x line=center,
        axis y line=center,
        scale only axis,
        clip bounding box=upper bound,
        clip=true,
        font=\footnotesize,
        restrict y to domain=-10:10,
    ]
    
    \addplot[color=color3, thick] 
        {1-(8 * ((2 * x^2 - 26) / 24)^4 - 8 * ((2 * x^2 - 26) / 24)^2 + 1) / 2.6300154321};
    \addplot[color=color6, thick] 
        {1-(8 * ((2 * x - 6) / 4)^4 - 8 * ((2 * x - 6) / 4)^2 + 1) / 23.5};
    \addplot[color=color3, thick, dashed] 
        {1.11022302e-16 + 0.0 *x + 1.15375705 *x^2 + 0.0 *x^3 - 0.34964182 *x^4 +
0.0 *x^5 + 0.03968907 *x^6 + 0.0 *x^7 - 0.00185947 *x^8 + 0.0 *x^9 +
(3.04831581e-05) *x^10};
    \addplot[color=color2, thick] 
        {1 - (1 - x^2 / 25)^4};

    \addplot[color=black, thick, dashed, domain=1:5] {1+1/23.5};
    \addplot[color=black, thick, dashed, domain=1:5] {1-1/23.5};
    \addplot[color=black, thick, dashed, domain=-5:-1] {1+1/2.6300154321};
    \addplot[color=black, thick, dashed, domain=-5:-1] {1-1/2.6300154321};
    \addplot[color=black, mark=*, only marks] coordinates {(0, 0)};

    \end{axis}
\end{tikzpicture}
\subcaption{
The interval $[-1, 1]$ is mapped to $[\lambda_{\min}, \lambda_{\max}]$. The polynomial is flipped and scaled to cross at $(0, 0)$.
}
\end{subfigure}
\hspace{0.03\textwidth}
\begin{subfigure}[c]{0.3\textwidth}
    \begin{tikzpicture}
    \begin{axis}[
        width=\textwidth,
        height=\textwidth,
        xlabel={$y$},
        domain=0:5,
        xtick=\empty,
        ytick=\empty,
        extra y ticks={0.8,1},
        extra x ticks={1,5},
        extra x tick labels={$\lambda_{\min}$, $\lambda_{\max}$},
        xmin=-1, xmax=5.5, 
        ymin=-0.36, ymax=2.2,
        samples=200,
        axis x line=center,
        axis y line=center,
        scale only axis,
        clip bounding box=upper bound,
        clip=true,
        font=\footnotesize,
        legend style={at={(1,0.95)},anchor=north east}
    ]
        \addplot[color=black, thick, domain=0.1:5] {1/x};
    \addlegendentry{$1/y$}
        \addplot[color=color2, thick, domain=-0.5:5] 
        {(1 - (1 - x^2 / 25)^4) / x};
    \addlegendentry{$P_{\qsvt}^{(7)}(y)$}
    \addplot[color=color6, thick, domain=-0.5:5] 
        {(1-(8 * ((2 * x - 6) / 4)^4 - 8 * ((2 * x - 6) / 4)^2 + 1) / 23.5) / x};
    \addlegendentry{$P_{\cheb}^{(3)}(y)$}
    \addplot[color=color3, thick,domain=-0.5:5] 
        {(1-(8 * ((2 * x^2 - 26) / 24)^4 - 8 * ((2 * x^2 - 26) / 24)^2 + 1) / 2.6300154321) / x};
    \addlegendentry{$P_{\qcheb}^{(7)}(y)$}
    \addplot[color=color3, thick, dashed, domain=-0.5:5] 
        {0.0 + 1.15375705 *x + 0.0 *x^2 - 0.34964182 *x^3 + 0.0 *x^4 +
0.03968907 *x^5 + 0.0 *x^6 - 0.00185947 *x^7 + 0.0 *x^8 +
(3.04831581e-05) *x^9};
    \addlegendentry{$P_{\chebopt}^{(7)}(y)$}

    \end{axis}
\end{tikzpicture}
\subcaption{
The final polynomial is obtained through division by $y$.\\
\label{fig:chebyshev-iteration-poly}
}
\end{subfigure}
    \caption{Constructing the polynomial of the Chebyshev iteration, the QSVT, $\qcheb$, and $\chebopt$ solvers for $n = 4$. The polynomial in (\textsc{b}) is $0$ at $y = 0$, and so division by  $y$ yields a polynomial. At the same time, it is close to $1$ for $y \in [\lambda_{\min}, \lambda_{\max}]$, meaning $1/y$ is well approximated in that interval.}
    \label{fig:chebyshev-iteration}
\end{figure}

It was a similar observation that led to the discovery of the \emph{Chebyshev iteration} for classical linear system algorithms~\cite{FS50}. Given bounds $0 < \lambda_{\min} < \lambda_{\max}$ on the eigenvalues of the matrix $A$, the Chebyshev iteration uses
\begin{equation} \label{eq:chebyshev-iteration}
P_{\cheb}^{(n)}(y) = \left(1 - \frac{T_{n+1}\bigl(\tfrac{2y - \lambda_{\min} - \lambda_{\max}}{\lambda_{\max} - \lambda_{\min}}\bigr)}{T_{n+1}\bigl(-\tfrac{\lambda_{\max} + \lambda_{\min}}{\lambda_{\max} - \lambda_{\min}}\bigr)}\right)/y,
\end{equation}
where again $T_n$ is the $n$-th Chebyshev polynomial; see~\eqref{eq:chebyshev}.
\Cref{fig:chebyshev-iteration} illustrates the underlying construction. For classical solvers, the choice~\eqref{eq:chebyshev-iteration} is known to be optimal with respect to~\eqref{eq:argmin} in the absence of any further information about $A$, see e.g.~\cite[Sec.~12.3.2]{Saa03}. It corresponds to the interpolation of $1/y$ with respect to the so-called \emph{Chebyshev nodes} in the interval $[\lambda_{\min}, \lambda_{\max}]$. Interpolation with respect to Chebyshev nodes is known to be particularly stable, which, together with~\eqref{eq:argmin}, motivates the choice heuristically.

In the quantum setting, the main obstacle in using the classical Chebyshev iteration polynomial is not its approximation quality on $[\lambda_{\min},\lambda_{\max}]$, but the fact that its magnitude need not remain small on the larger domain relevant for QSVT. This leads to increased normalization and hence higher cost.

To address this, \cite{GKS24} introduced a symmetrized variant of the Chebyshev iteration polynomial. The idea is to replace the dependence on $y$ inside the Chebyshev polynomial with $y^2$, while retaining the outer division by $y$. Concretely,
\[
P_{\qcheb}^{(2n-1)}(y) = \left(1 - {T_n\left(\frac{2y^2 - \lambda_{\min}^2 - \lambda_{\max}^2}{\lambda_{\max}^2 - \lambda_{\min}^2}\right)} / {T_n\left(-\frac{\lambda_{\max}^2 + \lambda_{\min}^2}{\lambda_{\max}^2 - \lambda_{\min}^2}\right)}\right) / y.
\]

Finally, \cite{SNW+25} gives an analytic expression for the polynomial minimizing the supremum-norm difference to $1/y$. Concretely, this polynomial results from replacing the Chebyshev polynomials in $P_\qcheb$ by
\[
L_n(x) = T_n(x) + \frac{1 - \lambda_{\min}}{1 + \lambda_{\min}} T_{n-1}(x),
\]
resulting in
\[
P_{\chebopt}^{(2n - 1)} = \left(1 - {L_n\left(\frac{2y^2 - \lambda_{\min}^2 - \lambda_{\max}^2}{\lambda_{\max}^2 - \lambda_{\min}^2}\right)} / {L_n\left(-\frac{\lambda_{\max}^2 + \lambda_{\min}^2}{\lambda_{\max}^2 - \lambda_{\min}^2}\right)}\right) / y.
\]

These constructions already account for normalization in the larger QSVT domain by virtue of their symmetric structure, but they do not optimize for it explicitly. It is not clear, however, whether this balance is optimal. Since the complexity bounds in \cref{lem:estimation,lem:qsvt} depend on both the approximation error and the normalization, it is natural to ask for polynomials that optimize this tradeoff directly. This motivates the constrained optimal polynomial constructions introduced next.

\section{Constrained uniform polynomial solver} \label{sec:cup}

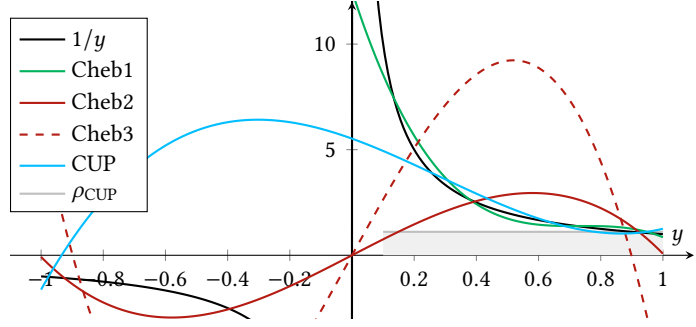
\begin{figure}
    \centering
    \begin{tikzpicture}
    \begin{axis}[
        width=0.6\textwidth,
        height=0.28\textwidth,
        xlabel={$y$},
        domain=-1:1,
        xmin=-1.1, xmax=1.1, 
        ymin=-3, ymax=12,
        samples=200,
        axis x line=center,
        axis y line=center,
        scale only axis,
        clip bounding box=upper bound,
        clip=true,
        font=\footnotesize,
        legend cell align=left,
        legend style={at={(0,0.95)},anchor=north west},
    ]

    \addplot[color=lightgray, thick,domain=0.1:1,samples=400,name path=rho, forget plot] {
1.0/0.9 * x + 1.0/0.9 * (1.0-x)
};

    \addplot[name path=floor, draw=none, forget plot] coordinates {(-1,0) (1,0)};
    
    \addplot[color=superlightgray, forget plot] fill between[of=rho and floor, soft clip={domain = 0.1:1}];

    \addplot[color=black, thick] {1/x};
    \addlegendentry{$1/y$}
    
    \addplot[color=color6, thick] {
1.73847751 - 2.29676906 *(-1.0 + 2.0*x) + 4.94688722 *(-1.0 + 2.0*x)^2 -
3.53349087 *(-1.0 + 2.0*x)^3
};
    \addlegendentry{$\cheb$}
    
    \addplot[color=color3, thick] {
2.86765399 + 0.98103952 *(-1.0 + 2.0*x) - 2.82992171 *(-1.0 + 2.0*x)^2 -
0.94330724 *(-1.0 + 2.0*x)^3
};
    \addlegendentry{$\qcheb$}

    \addplot[color=color3,dashed, thick] {
    9.21487603 + 0.95041322 *(-1.0 + 2.0*x) - 12.39669421 *(-1.0 + 2.0*x)^2 -
4.1322314 *(-1.0 + 2.0*x)^3
};
    \addlegendentry{$\chebopt$}

    \addplot[color=color4, thick] {
2.28567895 - 2.97715835 *(-1.0 + 2.0*x) + 1.10398266 *(-1.0 + 2.0*x)^2 +
0.84367226 *(-1.0 + 2.0*x)^3
};
    \addlegendentry{CUP}
        \addplot[color=lightgray, thick,domain=0.1:1,samples=400,name path=rho, draw=none, forget plot] {
1.0/0.9 * x + 1.0/0.9 * (1.0-x)
};
        
    \addplot[color=lightgray, thick,domain=0.1:1,samples=400,draw=none] {1};
    
    \addlegendentry{$\rho_\cups$}

    \end{axis}
\end{tikzpicture}
    \caption{
    Comparison of different approximations to $1/y$ with $\lambda_{\min} = 0.1$ and $\lambda_{\max} = 1$, as well as polynomial degree~$m = 3$. The Chebyshev iteration offers the best approximation, but is not constrained at all for negative values. The $\qcheb$ solver fixes this issue, but while the approximation for $y = 1$ is good, it is far from optimal for the smallest eigenvalue. The $\chebopt$ polynomial balances the approximation at $\lambda_{\min}$ and $\lambda_{\max}$, but still its approximation is hindered by its symmetry. The $\cups$ solver still offers a more constrained polynomial, while approximating $1/y$ better for positive $y$.}
    \label{fig:poly-comparison}
\end{figure}

We now turn the tradeoff between approximation quality and normalization into an explicit optimization problem. As in the previous constructions, we assume that spectral bounds~$\lambda_{\min}$ and~$\lambda_{\max}$ are available. Our starting point is the optimality criterion~\eqref{eq:argmin}. 
When the dimension~$d$ is large, it is natural to reinterpret the sum therein as an integral with respect to a spectral measure:
\begin{equation} \label{eq:rho1}
\sum_{j=1}^d (P(\lambda_j)-\lambda_j^{-1})^2 (V^\top b)_j^2
=
\int_{\lambda_{\min}}^{\mathrlap{\lambda_{\max}}} (P(y)-y^{-1})^2\,\mathrm{d}\rho(y),
\quad\text{where}\quad
\rho=\sum_{j=1}^d \delta(\bullet-\lambda_j)(V^\top b)_j^2.
\end{equation}
Here, $\delta$ denotes the Dirac measure. Since the measure~$\rho$ is not known explicitly, we replace it with a baseline model depending only on the available spectral interval, namely the uniform density
\[
\rho_\cups(y) = (\lambda_{\max}-\lambda_{\min})^{-1}
\qquad\text{for } y\in[\lambda_{\min},\lambda_{\max}].
\]

To obtain a tractable objective, we multiply the integrand in~\eqref{eq:rho1} by $y$, which cancels the $y^{-2}$~singularity:
\[
\int_{\lambda_{\min}}^{\lambda_{\max}} (P(y)-y^{-1})^2\, y\,\rho_\cups(y)\,\mathrm{d}y
=
\underbrace{\int_{\lambda_{\min}}^{\lambda_{\max}} y^{-1}\,\rho_\cups(y)\,\mathrm{d}y}_{\mathrm{const}}
+
\int_{\lambda_{\min}}^{\lambda_{\max}} \bigl(yP(y)^2-2P(y)\bigr)\,\rho_\cups(y)\,\mathrm{d}y.
\]
This is, up to the choice of spectral density, the same objective minimized by the Conjugate Gradient method in the classical setting; we return to this connection in~\cref{sec:cap}. In the quantum setting, however, approximation quality is only one part of the picture. The complexity bounds in \cref{lem:estimation,lem:qsvt} also depend on the normalization of the block encoding of the approximate solution. Assuming for simplicity that $\gamma_b=1$, this normalization is given by
\[
\gamma_x
=
\|P_{\mathrm{even}}\|_{[-\gamma_A,\gamma_A]}
+
\|P_{\mathrm{odd}}\|_{[-\gamma_A,\gamma_A]},
\]
see \cref{lem:qsvt}. Thus, the polynomial should approximate $1/y$ on the spectral interval~$[\lambda_{\min},\lambda_{\max}]$, while at the same time remaining small on the larger interval $[-\gamma_A,\gamma_A]$ relevant for QSVT; see \cref{fig:poly-comparison}.

To encode this tradeoff in a finite-dimensional optimization problem, we approximate the supremum norms of the even and odd parts by sampling on a uniform grid. More precisely, for every polynomial $Q$ with $\deg Q \le n$,
\[
\max_{-2n \le j \le 2n} \bigl|Q\bigl(\tfrac{j \gamma_A}{2n}\bigr)\bigr|
\le
\|Q\|_{[-\gamma_A, \gamma_A]}
\le
\tfrac43
\max_{-2n \le j \le 2n} \bigl|Q\bigl(\tfrac{j \gamma_A}{2n}\bigr)\bigr|.
\]
A similar method was considered in \cite{DLT22}, albeit with a smaller number of nodes and for minimizing the sup-norm error, not the normalization.
This allows us to enforce the sup-norm constraint via linear inequality constraints while keeping the block encoding normalization of the resulting polynomial controlled to within a factor of $4/3$.

We therefore compute $P_\mathrm{CUP} = P_{\rho_\cups}$ as the solution to the constrained minimization problem
\begin{multline}\label{eq:opt} 
P_{\rho} = \argmin_{\substack{\gamma_\mathrm{even},\, \gamma_\mathrm{odd} \in \R \\ \deg P \le n}} \tol(\gamma_\mathrm{odd} + \gamma_\mathrm{even})^2 + \sum_{j = 1}^{\Nquadrature} w_j \rho(y_j)(y_j P(y_j)^2 - 2P(y_j))\qquad \text{subject to} \\
-\gamma_\mathrm{even}\le P_\mathrm{even}\bigl(\tfrac{j \gamma_A}{2n}\bigr) \le \gamma_\mathrm{even},\quad
-\gamma_\mathrm{odd}\le P_\mathrm{odd}\bigl(\tfrac{j \gamma_A}{2n}\bigr) \le \gamma_\mathrm{odd}, 
\quad j = -2n, -2n+1, \ldots 2n.
\end{multline}
Here, $\tol>0$ is the target measurement accuracy, and $w_1,\dots,w_{\Nquadrature}$ and $y_1,\dots,y_{\Nquadrature}\in[\lambda_{\min},\lambda_{\max}]$ are the weights and nodes of Gauss--Legendre quadrature. Choosing $\Nquadrature = n+1$ suffices to evaluate the integral exactly. Polynomials obtained in this way will be referred to as \emph{Constrained Uniform Polynomial} (CUP) solvers.

The optimization problem~\eqref{eq:opt} is a convex constrained least-squares problem, and the unique solution can therefore be efficiently computed with standard numerical optimization tools. In our experiments, we use Sequential Least Squares Programming (SLSQP)~\cite{Kra88}.

Like the original QSVT solver, the $\qcheb$ and $\cups$ solvers fall into the class of \emph{semi-iterative} linear solvers~\cite[Section~5.5.3]{LS12}. Their accuracy improves with the number of steps, but the approximation interval $[\lambda_{\min},\lambda_{\max}]$ must be specified a priori.

\section{Constrained adaptive polynomial solver} \label{sec:cap}
In classical numerical linear algebra, the fastest polynomial solvers are typically not fixed semi-iterative methods but adaptive Krylov methods such as the conjugate gradient (CG). Rather than relying on a priori spectral bounds $\lambda_{\min}$ and $\lambda_{\max}$, these methods extract problem-specific spectral information during the iteration and use it to build improved approximation polynomials.

We aim to transfer this idea to the quantum setting by introducing a \emph{Constrained Adaptive Polynomial} (CAP) quantum linear solver. At a high level, CAP first estimates quantities of interest related to the linear system, namely Chebyshev moments. It then finds a distribution of eigenvalues consistent with these moments. Finally, it computes an optimal polynomial given this distribution while taking into account the effect of normalization on measurement accuracy.

\begin{figure}
    \centering
    \begin{tikzpicture}
    \begin{axis}[
        width=0.6\textwidth,
        height=0.28\textwidth,
        xlabel={$y$},
        domain=-1:1,
        xmin=-1.1, xmax=1.1, 
        ymin=-1.5, ymax=12,
        samples=200,
        axis x line=center,
        axis y line=center,
        scale only axis,
        clip bounding box=upper bound,
        clip=true,
        font=\footnotesize,
        legend cell align=left,
        legend style={at={(0,0.95)},anchor=north west},
    ]

    \addplot[color=lightgray, thick,domain=0.1:1,samples=400, name path=rho, forget plot] {
exp(
-4.19786972 + 1.49212554 *(-1.0 + 2.0*x) + 60.11959054 *(-1.0 + 2.0*x)^2 +
35.60711375 *(-1.0 + 2.0*x)^3 - 150.74350797 *(-1.0 + 2.0*x)^4 -
169.41801626 *(-1.0 + 2.0*x)^5 - 61.81019411 *(-1.0 + 2.0*x)^6 -
7.54908737 *(-1.0 + 2.0*x)^7
)
};

    \addplot[name path=floor, draw=none, forget plot] coordinates {(-1,0) (1,0)};
    
    \addplot[color=superlightgray, forget plot] fill between[of=rho and floor, soft clip={domain = 0.1:1}];

    \addplot+[
        scatter,
        color=black,
        mark=o,
        only marks,
        scatter src=explicit,
        scatter/use mapped color={fill=gray,draw=black},
        scatter/@pre marker code/.append style=
                {/tikz/mark size={.6pt+\pgfplotspointmetatransformed/650}}
    ] table[x=x, y=y, meta=z, row sep=crcr] {
        x y z \\
        2.462695355026068011e-01 4.060591570772727366e+00 -2.129363065249910945e-03\\
7.061817494706966736e-01 1.416066049213999811e+00 1.297474292503504600e-01\\
2.789746676434467498e-01 3.584554857425566343e+00 -8.193109472914517077e-02\\
6.965959571413917750e-01 1.435552402720914333e+00 -1.815522685738422959e-01\\
2.646682623204057450e-01 3.778314752334778426e+00 -3.739465436801709397e-01\\
6.682150967457112234e-01 1.496524105591330756e+00 -5.538193326800023852e-01\\
6.692454126143490001e-01 1.494220178654026121e+00 2.240174054323143260e-01\\
7.020570200424232077e-01 1.424385728583090094e+00 -1.678444782196305829e-01\\
7.463176549388931269e-01 1.339912024567980753e+00 2.079761064479696553e-01\\
6.394942944488142267e-01 1.563735609028238160e+00 7.054942605130834488e-02\\
7.888188762419912070e-01 1.267718141792062436e+00 -1.489701145615239808e-01\\
2.224729391925784339e-01 4.494928702921363595e+00 -2.129895288578766421e-01\\
7.507303519614689780e-01 1.332036193004921598e+00 7.102493198723418277e-02\\
2.496097028062316814e-01 4.006254519586064156e+00 -1.070008889878617392e-01\\
2.676574252336847426e-01 3.736119030237722516e+00 -9.541452545300248278e-02\\
7.079224721720974722e-01 1.412584060132643815e+00 -3.141351348756466627e-02\\
6.885012680789824868e-01 1.452430149896664302e+00 1.605313535377585199e-01\\
7.581400909245967767e-01 1.319017437503457568e+00 -2.304369425628114698e-01\\
6.967560382101620764e-01 1.435222581735805081e+00 -3.862548950265919934e-01\\
7.095353991808197591e-01 1.409372951870379476e+00 2.194150630966749949e-01\\
    };
    \addlegendentry{$\lambda_0, \dots, \lambda_{d-1}$}
    \addplot[color=black, thick] {1/x};
    \addlegendentry{$1/y$}

    \addplot[color=color5, thick] {
1.8523302 - 1.65833324 *(-1.0 + 2.0*x) + 3.34408457 *(-1.0 + 2.0*x)^2 -
4.45575441 *(-1.0 + 2.0*x)^3
};
    \addlegendentry{CG}

        \addplot[color=color4, thick] {
2.28567895 - 2.97715835 *(-1.0 + 2.0*x) + 1.10398266 *(-1.0 + 2.0*x)^2 +
0.84367226 *(-1.0 + 2.0*x)^3
};
    \addlegendentry{CUP}

        \addplot[color=color1, thick] {
2.59680649 - 2.43704603 *(-1.0 + 2.0*x) - 1.21852302 *(-1.0 + 2.0*x)^2
};
    \addlegendentry{$\caps$}
    
    \addplot[color=lightgray, thick,domain=0.1:1,samples=400,draw=none] {1};
    \addlegendentry{$\rho_\caps$}
    \end{axis}
\end{tikzpicture}
    \caption{Comparison of different approximations to $1/y$ with $\lambda_{\min} = 0.1$ and $\lambda_{\max} = 1$, as well as polynomial degree~$m = 3$. The size of the circles indicates the weight $(V^\top b)_j$ of the eigenvalues as given in~\eqref{eq:argmin}. The CG polynomial well approximates $1/y$ at the clusters of eigenvalues, but increases quickly for negative values. The $\cups$ solver has no knowledge of the distribution of eigenvalues and aims at uniform approximation on the interval~$[0.1, 1]$ while restricting its global maximum. The $\caps$ solver recovers a distribution~$\rho_\caps$ of eigenvalues consistent with the moments and concentrates its approximation of $1/y$ on the eigenvalue clusters, at the cost of accuracy elsewhere, which is exactly the tradeoff CG exploits classically.}
    \label{fig:poly-comparison2}
\end{figure}
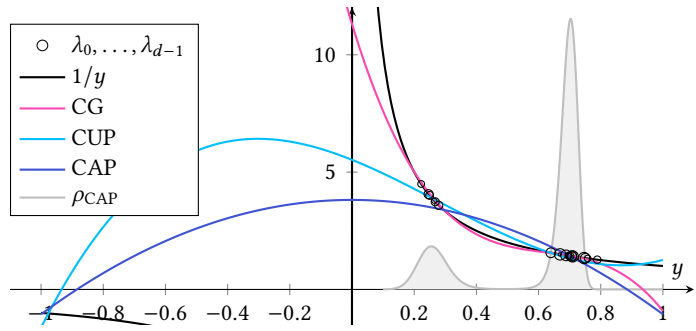

In the limit of exact moment measurements and without the normalization constraint, the polynomial obtained through this process coincides with the one generated implicitly by classical conjugate gradients (CG); see \cite[Sec.~6.7]{Saa03}. 
Classically, CG can exploit this polynomial through short recurrences, without explicitly computing spectral information. On a quantum computer, however, this recurrence cannot be exploited in the same way: each inner product requires a separate state preparation, since the No-Cloning Theorem~\cite{WZ82} forbids reusing an already-measured state. An explicit spectral reconstruction is therefore natural in the quantum setting.
The central design decision is then how to use this spectral information. A polynomial adapted purely to the spectrum, as in CG, can grow rapidly outside the spectral interval and thereby inflate the block encoding normalization, which directly increases measurement cost via~\cref{lem:estimation} (see the purple curve in~\cref{fig:poly-comparison2}). The $\cups$ solver of the previous section already addresses this issue by coupling approximation accuracy with a supremum-norm constraint on the polynomial. The $\caps$ solver adopts the same constrained optimization but replaces CUP's uniform spectral model with the reconstructed density~$\rho_\caps$, making the solver adaptive to the actual spectrum of the problem while preserving CUP's control on normalization.

\subsection{Moment estimation} \label{ssec:bending}
The first step in the $\caps$ solver is to extract spectral information from the linear system. For this, we estimate quantities of the form
\[\mu_k = b^\top R_k(A) b \quad\text{for}\quad k = 0, \dots, 2n+1,\]
which we call \emph{moments}, and where $R_0, \dots, R_{2n+1}$ form a basis of the polynomials of degree at most~\mbox{$2n+1$}. Once these moments are known, Euclidean and $A$-inner products between vectors of the form~$P(A)b$ could be computed classically, without further quantum measurements.

The simplest choice is to take $R_k(y) = y^k$, so that $\mu_k = b^\top A^k b$. We instead use \emph{Chebyshev moments}, with $R_k$ given by rescaled Chebyshev polynomials
\[
R_k(y) = T_k(y / \gamma_A).
\]
These are efficiently computable via QSVT~\cite[Lemma~9]{GSLW19} and yield a better-conditioned characterization of the spectral information relevant to the solver. In our experiments, Chebyshev moments performed near-optimally among the polynomial bases considered.\footnote{This was studied using genetic algorithms to find optimal polynomials. The results are available at \url{https://github.com/MDeiml/quantum-krylov/blob/main/test_basis.ipynb}.}

The final $\caps$ polynomial has degree $n$, corresponding to a circuit of depth $\calO(n\, T_A)$. Moments up to degree $2n+1$ can be obtained within the same depth budget. Using the technique of~\cite{KMM23}, the $k$-th Chebyshev moment is computed in $\lceil k/2 \rceil$ applications of the block encoding of $A$. This relies on the fact that the QSVT phase angles for Chebyshev polynomials are symmetric. The same technique extends to any quantity of the form $b^\top P(A) b$ whenever symmetric phase angles for $P$ can be computed, e.g., via~\cite{DMWL21,NSYL25}.

\subsection{Computing the maximum entropy distribution} \label{ssec:supnorm}

Given the moments $\mu_0, \dots, \mu_{2n + 1}$, there are many measures that exactly match this data. The simplest choice is to model $\rho$ as a weighted sum of Dirac measures, as in \eqref{eq:rho1}, which would yield the so-called \emph{Ritz values}. We found, however, that recovering Ritz values from Chebyshev moments is numerically unstable, as the moment-to-nodes problem is ill-conditioned even under moderate measurement noise. Standard stabilization procedures, such as~\cite[Algorithm~1]{ZWXL24} or~\cite[Algorithm~1.1]{ELN22}, trade this instability for a loss of spectral coverage: components of the Krylov subspace carrying non-negligible weight are discarded, degrading overall solver performance.

For this reason, we chose to model $\rho$ as a continuous measure, namely the maximum entropy distribution~\cite{MP84} with the given moments. Its density~$\rho_\caps \colon [\lambda_{\min}, \lambda_{\max}] \to [0, \infty)$ is given by
\begin{align*}
\rho_\caps \coloneq\ &\argmax_{\rho \ge 0} \, -\int_{\lambda_{\min}}^{\lambda_{\max}} \rho(y) \log \rho(y) \, \mathrm{d}y \\
&\text{subject to}\quad
\int_{\lambda_{\min}}^{\lambda_{\max}} y^j \rho(y) \, \mathrm{d}y = \mu_j \quad\text{for } j = 0, \dots, 2n+1.
\end{align*}
The corresponding measure $\mathrm{d}\rho_\caps(y) = \rho_\caps(y)\,\mathrm{d}y$ can then be substituted for the discrete spectral measure in~\eqref{eq:rho1}.
The max-entropy problem is well-posed on any interval containing the spectrum, e.g., $[0, \gamma_A]$ or $[-\gamma_A, \gamma_A]$ in the indefinite case. Restricting to $[\lambda_{\min}, \lambda_{\max}]$ incorporates the available spectral bounds and yields a more concentrated reconstruction.
A standard Lagrangian argument shows that the solution has the form
\[
\rho_\caps(y) = \exp(Q(y))
\]
with polynomial~$Q$ of degree at most $2n+1$. It remains to find the coefficients of $Q$ for which the moments are reproduced~\cite{MP84}. In our numerical experiments, we use the BFGS method for this purpose. We again use Gauss--Legendre quadrature as in~\eqref{eq:opt} to approximate the integral. However, the number of quadrature points must be increased beyond $\Nquadrature = n + 1$, since the integrand is no longer a polynomial. In our experiments, we use $\Nquadrature = 4(n+1)$, which was sufficient for all problems considered.
Using~\eqref{eq:opt} to compute the optimal polynomial $P_\caps = P_{\rho_\caps}$ with the resulting density completes the $\caps$ solver.

\section{Structural transformations of solver polynomials} \label{sec:transforms}

The polynomial constructions introduced so far can be combined with additional transformations that exploit structure in the matrix representation. These transformations are largely independent of the specific choice of solver polynomial and may be used, for example, with the QSVT, $\qcheb$, CUP, or $\caps$ polynomials. Given a solver polynomial $P$, we denote by $P_\mathrm{sq}$ and $P_\mathrm{sq}'$ its inner and outer square transforms. The inner transform squares the argument inside $P$, while the outer transform squares $P$ itself.

Both transformations discussed below assume access to a block encoding of a (potentially non-square) matrix $B$ such that
\[
A = B^\top B.
\]
The matrix $B$ need not be symmetric or even square. Indeed, QSVT on $B$ with an even polynomial behaves as if $B = A^{1/2}$, which we will assume here for simplicity. While this is not immediate from \cref{lem:qsvt}, it follows from the general QSVT framework; see~\cite[Section~3.2]{GSLW19}. Whether such a~$B$ admits an efficient block encoding depends on the problem structure. In many settings, it arises naturally. For instance, from a Cholesky-like decomposition or from the underlying discretization of a differential operator; see~\cite{DP25} for a concrete example.

For the first transformation, we consider the transformed polynomial
\begin{equation} \label{eq:square-transform}
P_\mathrm{sq}(B) \coloneq P(B^2) = P(A).
\end{equation}
This transformation replaces the normalization domain $[-\gamma_B,\gamma_B]$ with $[0,\gamma_B^2]$ since
\[
\|P_\mathrm{sq}\|_{[-\gamma_B,\gamma_B]} = \|P\|_{[0,\gamma_B^2]}.
\]
In particular, the negative half-axis no longer contributes to the normalization, while the behavior on $[0,\lambda_{\min})$ remains unchanged.

A second transformation, also considered in~\cite{DP25}, uses the identity
\[
A^{-1} = (B^{-1})^2.
\]
As any eigenvalue~$\lambda_B$ of $B$ satisfies $\lambda_{\min} \le \lambda_B^2 \le \lambda_{\max}$, the polynomial~$P$ only needs to approximate~$1/y$ on this smaller spectrum.
If $P$ is additionally of odd parity, that is, if it approximates $1/y$ also for negative values, one may define
\begin{equation} \label{eq:square-outer-transform}
P_\mathrm{sq}'(B) \coloneq P(B)^2 \approx (B^{-1})^2 = A^{-1}.
\end{equation}
This transformation also changes the behavior near $y=0$, since odd parity implies $P_\mathrm{sq}'(0)=0$. 

A third natural alternative is to shift the matrix so that its spectrum becomes symmetric. This idea is considered in~\cite{TWYH24}, where one replaces $A$ with $A - \lambda \Id$ using $\lambda \coloneq (\lambda_{\max} + \lambda_{\min}) / 2$
and modifies the polynomial accordingly. In general, however, such a shift does not improve efficiency by itself. Indeed, applying a polynomial $\tilde P$ to the shifted matrix $A - \lambda \Id$ is equivalent to applying the polynomial $y \mapsto \tilde P(y-\lambda)$ to $A$, so this transformation does not go beyond \cref{lem:qsvt}. Still, the shift may be useful in combination with additional techniques such as uniform singular value amplification~\cite[Theorem~30]{GSLW19}, or in sparse settings with large diagonal entries; see~\cite{TWYH24}.

\section{Numerical experiments}\label{sec:num}
We benchmark the proposed $\cups$ and $\caps$ solvers against the QSVT, $\qcheb$, and $\chebopt$ solvers using numerical simulations of the QSVT and estimation primitives. Code is available at \url{https://github.com/MDeiml/quantum-krylov}.

\subsection{Benchmark setup}
Our goal is to compare the solvers at the level of the abstract primitives used in their analysis, namely the estimation and QSVT procedures from \cref{lem:estimation,lem:qsvt}. We therefore use a custom simulation framework that models these primitives directly, rather than explicit quantum circuits for particular block encodings. This abstraction lets us vary the noise level parametrically and study performance across regimes ranging from near-term noisy devices to low-error fault-tolerant computation.

\paragraph{Benchmark suite}
We consider linear systems of the form $Dx = b$, where $b$ is a uniformly random vector with $\|b\|_2 = 1$, and $D = \diag(\lambda_0, \dots, \lambda_{d-1})$ is a random diagonal matrix. For the distribution of its eigenvalues, we consider two cases. In both cases, the distribution is parameterized by an upper bound $\tilde \kappa$ on the resulting condition number; the realized $\kappa
$ is a random variable with $\kappa \le \tilde\kappa$.

In case $1$, we set $\lambda_0 = 1/\tilde{\kappa}$ and draw the other eigenvalues $\lambda_1, \dots, \lambda_{d-1}$ independently and uniformly from $[\lambda_0, 1]$. In case $2$, to model clustered eigenvalues, we pick $4$ cluster locations $\tilde \lambda_1, \dots, \tilde \lambda_4$ uniformly in $[1/\tilde{\kappa}, 1]$. Each eigenvalue is subsequently sampled from a normal distribution around a random cluster location with a standard deviation of $0.025 / (1 - 1/\tilde{\kappa})$. The resulting values are clamped to $[1/\tilde{\kappa}, 1]$. We reiterate that $\tilde{\kappa}$ refers to an upper bound, given as a parameter to the solver; the realized condition number $\kappa$ is typically slightly smaller.

Restricting to diagonal matrices does not reduce generality. A general SPD matrix with either of the above eigenvalue distributions has the form $A = UDU^\top$ with $U$ Haar-distributed, and QSVT applied to $A$ is unitarily equivalent to QSVT applied to $D$ with $b$ replaced by $U^\top b$. Since $U$ is Haar and $b$ is uniform on the unit sphere, $U^\top b$ has the same distribution as $b$, so the diagonal simulation samples the same ensemble of problems as the general one.

\paragraph{Noise model}
We simulate noise by randomly inserting Pauli flips in the simulation of QSVT. A Pauli flip in an orthonormal basis $U \in \C^{2^q \times 2^q}$ of the state space with $q \in \N$ qubits is given by
\[
F(U) \coloneq U \left(\begin{bmatrix}
    1 & 0 \\ 0 & -1
\end{bmatrix} \otimes \Id_{2^{q-1}} \right) U^\dag.
\]
Due to the randomness of the basis, the choice of flip -- $X$, $Y$, or $Z$ -- does not matter.
To simulate the matrix being block encoded, we choose a number of Haar-distributed unitaries~$U_1, \dots, U_{\Nnoise}$, representing the points in the circuit at which noise might occur. The noise level is controlled by a parameter $\xi \ge 0$, which determines the average number of noise events per application of the block encoding. In the simulation, after each application of the block encoding of $A$ and with probability $\xi$, a noise flip~$F(U)$, with $U$ randomly chosen from $U_1, \dots, U_{\Nnoise}$, is applied to the state vector.

\paragraph{Parameters}
We set $d = 128$, $\tilde{\kappa} = 3$, and $\Nnoise = 20$ and test all combinations of parameters where
\[
n \in \{1, 2, \ldots, 16\} ,\quad \samples \in \{1, 4, 16\} \cdot 10^4,\quad\text{and}\quad \xi \in \{0, 0.0025, 0.005, \dots, 0.04\}.
\]
These values are chosen to make the simulations tractable while keeping the problem non-trivial.
The $\qsvt$, $\cups$, and $\caps$ solvers each take a target accuracy $\tol$ as input, which we set as $\tol = 2/\sqrt{\samples}$, reflecting the Monte Carlo standard deviation of the estimation procedure from \cref{lem:estimation}. For $\qsvt$ this determines the reference parameter $\tilde n = \lceil \tilde{\kappa}^2 \log(\tilde{\kappa}/\tol) \rceil$ as suggested in~\cite[Lemma~40]{GSLW19}. For $\cups$ and $\caps$ it sets the weight of the normalization penalty in~\eqref{eq:opt}. The $\qcheb$ and $\chebopt$ solvers do not take $\tol$ as input, as their degree is chosen directly via $n$.

\subsection{Benchmark results}
\begin{figure}
\begin{tikzpicture}
    \begin{groupplot}[group style={group size=3 by 2,
      horizontal sep=0pt,vertical sep=0pt,xticklabels at=edge bottom},
        ymax = 0.8,
        ymin = 0.003,
        ymode=log,
        height=0.35\textwidth,
        xmin=0.1,xmax=15.5,
        width=0.38\textwidth,
        xtick={1,3,5,7,9,11,13,15},
    ]
    \nextgroupplot[
        font=\footnotesize,
        ylabel={$\err(A, b)$ uniform $\lambda$},
    ]
    \errorplot[color=color2,mark=square, solid]{qsvt_None}{data/uniform/0_160000.csv}{0.2}
    \errorplot[color=color3,mark=triangle, solid]{q_cheb_None}{data/uniform/0_160000.csv}{0}
    \errorplot[color=color3,mark=pentagon, dashed]{chebopt_None}{data/uniform/0_160000.csv}{0}
    \errorplot[color=color4,mark=x]{cap_None_False}{data/uniform/0_160000.csv}{-0.2}
    \errorplot[color=color1,mark=o]{cap_None_True}{data/uniform/0_160000.csv}{-0.2}
    \nextgroupplot[
        font=\footnotesize,
        yticklabels=\empty,
    ]
    \errorplot[color=color2,mark=square, solid]{qsvt_None}{data/uniform/0.0025_160000.csv}{0.2}
    \errorplot[color=color3,mark=triangle, solid]{q_cheb_None}{data/uniform/0.0025_160000.csv}{0}
    \errorplot[color=color3,mark=pentagon, dashed]{chebopt_None}{data/uniform/0.0025_160000.csv}{0}
    \errorplot[color=color4,mark=x]{cap_None_False}{data/uniform/0.0025_160000.csv}{-0.2}
    \errorplot[color=color1,mark=o]{cap_None_True}{data/uniform/0.0025_160000.csv}{-0.2}
    \nextgroupplot[
        legend style={
            at={(1.0,1.02)},anchor=south east,
        },
        legend columns=5,
        font=\footnotesize,
        yticklabel pos=right,
    ]
    \errorplot[color=color2,mark=square, solid]{qsvt_None}{data/uniform/0.005_160000.csv}{0.2}
    \errorplot[color=color3,mark=triangle, solid]{q_cheb_None}{data/uniform/0.005_160000.csv}{0}
    \errorplot[color=color3,mark=pentagon, dashed]{chebopt_None}{data/uniform/0.005_160000.csv}{0}
    \errorplot[color=color4,mark=x]{cap_None_False}{data/uniform/0.005_160000.csv}{-0.2}
    \errorplot[color=color1,mark=o]{cap_None_True}{data/uniform/0.005_160000.csv}{-0.2}

    \addlegendentry{$\qsvt$}
    \addlegendentry{$\qcheb$}
    \addlegendentry{$\chebopt$}
    \addlegendentry{$\cups$}
    \addlegendentry{$\caps$}
    
    \nextgroupplot[
        font=\footnotesize,
        ylabel={$\err(A, b)$ clustered $\lambda$},
    ]
    \errorplot[color=color2,mark=square, solid]{qsvt_None}{data/clustered/0_160000.csv}{0.2}
    \errorplot[color=color3,mark=triangle, solid]{q_cheb_None}{data/clustered/0_160000.csv}{0}
    \errorplot[color=color3,mark=pentagon, dashed]{chebopt_None}{data/clustered/0_160000.csv}{0}
    \errorplot[color=color4,mark=x]{cap_None_False}{data/clustered/0_160000.csv}{-0.2}
    \errorplot[color=color1,mark=o]{cap_None_True}{data/clustered/0_160000.csv}{-0.2}
    \nextgroupplot[
        font=\footnotesize,
        xlabel={polynomial degree~$m$},
        yticklabels=\empty,
    ]
    \errorplot[color=color2,mark=square, solid]{qsvt_None}{data/clustered/0.0025_160000.csv}{0.2}
    \errorplot[color=color3,mark=triangle, solid]{q_cheb_None}{data/clustered/0.0025_160000.csv}{0}
    \errorplot[color=color3,mark=pentagon, dashed]{chebopt_None}{data/clustered/0.0025_160000.csv}{0}
    \errorplot[color=color4,mark=x]{cap_None_False}{data/clustered/0.0025_160000.csv}{-0.2}
    \errorplot[color=color1,mark=o]{cap_None_True}{data/clustered/0.0025_160000.csv}{-0.2}
    \nextgroupplot[
        legend style={
            at={(1.0,1.02)},anchor=south east,
        },
        legend columns=4,
        font=\footnotesize,
        yticklabel pos=right,
    ]
    \errorplot[color=color2,mark=square, solid]{qsvt_None}{data/clustered/0.005_160000.csv}{0.2}
    \errorplot[color=color3,mark=triangle, solid]{q_cheb_None}{data/clustered/0.005_160000.csv}{0}
    \errorplot[color=color3,mark=pentagon, dashed]{chebopt_None}{data/clustered/0.005_160000.csv}{0}
    \errorplot[color=color4,mark=x]{cap_None_False}{data/clustered/0.005_160000.csv}{-0.2}
    \errorplot[color=color1,mark=o]{cap_None_True}{data/clustered/0.005_160000.csv}{-0.2}

    \end{groupplot}
    \draw (group c1r1.north west) node[draw, anchor=north west, outer sep=0.2cm] {$\xi = 0$};
    \draw (group c2r1.north west) node[draw, anchor=north west, outer sep=0.2cm] {$\xi = 0.0025$};
    \draw (group c3r1.north west) node[draw, anchor=north west, outer sep=0.2cm] {$\xi = 0.005$};
\end{tikzpicture}
\caption{
Comparison of the $\qsvt$, $\qcheb$, $\chebopt$, $\cups$, and $\caps$ solvers with  $\xi = 0$ (left), $\xi = 0.0025$ (middle), and $\xi = 0.005$ (right), as well as $\samples = 16 \cdot 10^4$. The reported measure is the error from \cref{def:complexity} relative to the norm of the true solution. The dotted line indicates the 95th percentile of observed errors.
\label{fig:results}
}
\end{figure}
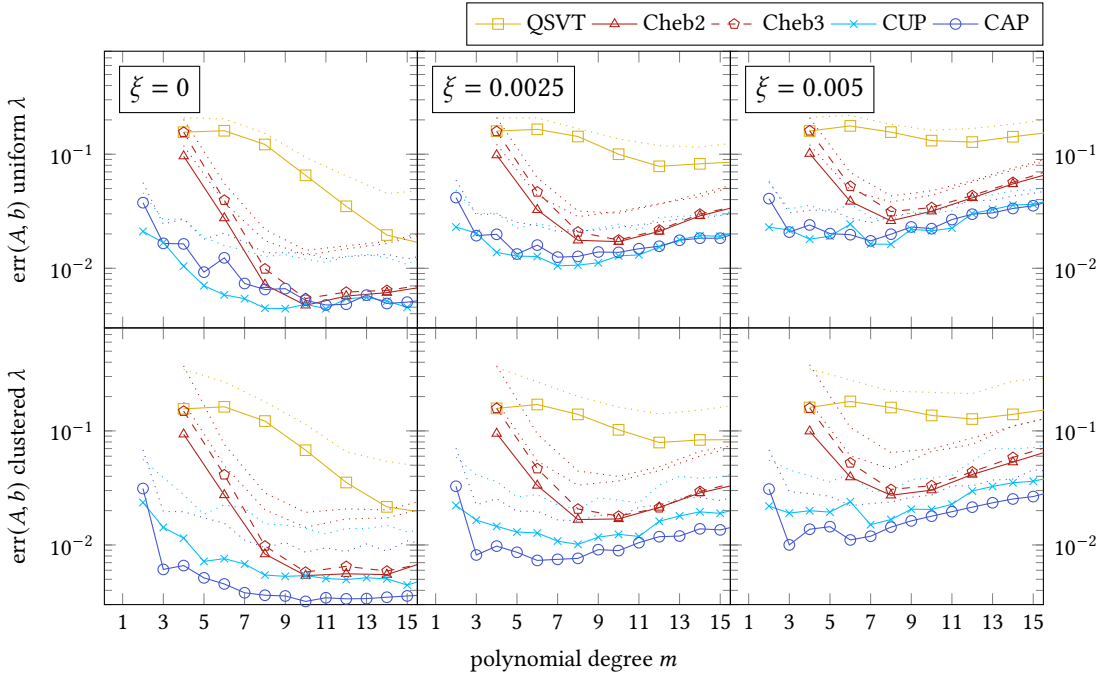

We simulate the $\qsvt$, $\cheb$, $\qcheb$~\cite{GKS24}, $\chebopt$~\cite{SNW+25}, $\cups$, and $\caps$ solvers described in \cref{sec:qsvt,sec:chebsym,sec:cup,sec:cap}, each with and without the inner and outer square transforms~\eqref{eq:square-transform} and~\eqref{eq:square-outer-transform}, denoted by subscripts $\ast_\mathrm{sq}$ and~$\smash{\ast_\mathrm{sq}'}$, respectively. For each parameter combination described above, we draw $200$ random test equations and report the median of the error~$\err(A, b)$ from \cref{def:complexity}. The $95^{\mathrm{th}}$ percentiles are shown as dotted lines in \cref{fig:results}. For each solver and noise level, the optimal polynomial degree~$n$ is selected retrospectively to minimize median error; convergence behavior as a function of~$n$ is shown separately in \cref{fig:results}. The full data set is available at \cite{DP26a}.

\begin{table}
    \centering
    \small
    \begin{tabular}{lllllllllll}
        \toprule
        &\multicolumn{4}{c}{uniform eigenvalues} &\multicolumn{4}{c}{clustered eigenvalues}\\ \cmidrule(l){2-5} \cmidrule(l){6-9}
        $\xi = $ & $0$ & $0.005$ & $0.01$ & $0.02$ & $0$ & $0.005$ & $0.01$ & $0.02$ \\
        \midrule
$\qsvt$ & 0.012(11) & 0.128(5) & 0.161(1) & 0.171(1) & 0.012(14) & 0.127(5) & 0.164(1) & 0.176(1) \\
$\qcheb$ & 0.005(4) & 0.026(3) & 0.047(3) & 0.073(2) & 0.005(4) & 0.027(3) & 0.045(3) & 0.074(2) \\
$\chebopt$ & 0.005(4) & 0.031(3) & 0.052(3) & 0.089(2) & 0.006(4) & 0.031(3) & 0.050(3) & 0.089(2) \\
$\cups$ & 0.004(10) & 0.016(7) & 0.021(1) & \textbf{0.024(1)} & 0.004(14) & 0.015(6) & 0.023(1) & 0.026(1) \\
$\caps$ & 0.005(10) & 0.017(6) & 0.027(2) & 0.036(2) & 0.003(9) & \textbf{0.010(2)} & 0.017(2) & 0.028(2) \\
\midrule
        $\xi = $ & $0$ & $0.0025$ & $0.005$ & $0.01$ & $0$ & $0.0025$ & $0.005$ & $0.01$ \\
        \midrule
$\qcheb_\mathrm{sq}'$ & 0.005(2) & 0.019(2) & 0.036(2) & 0.060(1) & 0.006(2) & 0.020(2) & 0.035(2) & 0.060(1) \\
$\chebopt_\mathrm{sq}'$ & 0.006(2) & 0.021(2) & 0.042(2) & 0.076(2) & 0.007(2) & 0.021(2) & 0.042(2) & 0.076(2) \\
$\cups_\mathrm{sq}$ & 0.004(10) & 0.015(1) & \textbf{0.016(1)} & 0.025(1) & 0.004(15) & 0.014(1) & 0.018(1) & 0.025(1) \\
$\caps_\mathrm{sq}$ & \textbf{0.004(16)} & \textbf{0.012(1)} & 0.018(1) & 0.026(1) & \textbf{0.003(9)} & 0.011(1) & \textbf{0.014(1)} & \textbf{0.023(1)} \\
        \bottomrule
\end{tabular}
    \caption{
    Best results for selected solvers, with $\samples = 16 \cdot 10^4$ and uniform as well as clustered eigenvalues. The results for the transformed solvers assume half of the noise rate for the square root~$B$, as would typically be the case, see e.g.\ \cite{DP25}. The reported quantity is the median of $\err(A, b)$ over all test equations. For each combination of solver and noise, the best number of steps~$n$ is selected, indicated in parenthesis after the result. The best solver in each column is highlighted.
    }
    \label{tab:best-results}
\end{table}

\Cref{tab:all-results,tab:all-results-clustered} contain a comparison of the minimal median error achieved using each method. All methods perform well without noise and with a large number of samples, i.e., $\xi = 0$ and $\samples = 16 \cdot 10^4$. For this noiseless case, decreasing the number of samples increases the error. As~$\xi$ increases, the error increases, and the differences between different values of $\samples$ disappear. 
\Cref{tab:best-results} contains a concise overview of the results with the largest sampling accuracy $\samples = 16 \cdot 10^4$. Transformed variants are shown only when they outperform the untransformed solver in at least one configuration. The performance of the different solvers is then as follows.
\begin{itemize}
    \item The $\qsvt$ solver performs reasonably well at $\xi = 0$, but it degrades substantially as noise increases.

    \item The $\cheb$ solver, which is not listed in \cref{tab:best-results}, performs surprisingly well with $n = 1$, but the error always increases for $n\leq 2$. While the transform \eqref{eq:square-outer-transform} cannot be applied here, the transform \eqref{eq:square-transform} improves performance, especially for small noise.

    \item The $\qcheb$ and $\chebopt$ solver have the best convergence behavior out of the analytically constructed methods. If a block encoding of a matrix square root~$B$ is available and has half the gate count of $A$, the transform~\eqref{eq:square-outer-transform} seems to improve their performance even further. $\qcheb$ seems to have slightly better minimal error for the test cases considered here. The gap increases for increasing noise, indicating that the theoretical optimality of $\chebopt$ holds only in the context of accurate measurements. Additional results for $\tilde \kappa = 5$ are available in the full dataset~\cite{DP26a} verify, that for $N = 16 \cdot 10^4$ and $\xi = 0$, the $\chebopt$ polynomial is indeed slightly better, as predicted by the results regarding sup-norm in~\cite{SNW+25}.

    \item The $\cups$ solver achieves the lowest median error in nearly all regimes with uniformly random eigenvalues. If a matrix square root is available, the transformed variant $\cups_\mathrm{sq}$ provides further improvement.
    \item The $\caps$ solver performs best regarding the median error in almost all regimes for clustered random eigenvalues. Again, the variant with~\eqref{eq:square-transform} should be used if possible.
\end{itemize}

\begin{figure}
\begin{tikzpicture}
    \begin{groupplot}[group style={group size=3 by 2,
      horizontal sep=0pt,vertical sep=0pt,xticklabels at=edge bottom,},
        ymax = 0.8,
        ymin = 0.003,
        xmin = 1.1e5,
        xmax = 0.9e8,
        xmode=log,
        ymode=log,
        height=0.35\textwidth,
        width=0.38\textwidth,
    ]
    \nextgroupplot[
        ylabel={$\err(A, b)$ uniform $\lambda$},
        font=\footnotesize,
    ]

    \complexityplot[color=color2,mark=square, solid]{qsvt_None}{data/uniform/0_160000.csv}
    \complexityplot[color=color3,mark=triangle, solid]{q_cheb_None}{data/uniform/0_160000.csv}
    \complexityplot[color=color3,mark=pentagon, dashed]{chebopt_None}{data/uniform/0_160000.csv}
    \complexityplot[color=color4,mark=x]{cap_None_False}{data/uniform/0_160000.csv}
    \complexityplot[color=color1,mark=o]{cap_None_True}{data/uniform/0_160000.csv}
    
    \nextgroupplot[
        font=\footnotesize,
        yticklabels=\empty,
    ]
    
    \complexityplot[color=color2,mark=square, solid]{qsvt_None}{data/uniform/0.0025_160000.csv}
    \complexityplot[color=color3,mark=triangle, solid]{q_cheb_None}{data/uniform/0.0025_160000.csv}
    \complexityplot[color=color3,mark=pentagon, dashed]{chebopt_None}{data/uniform/0.0025_160000.csv}
    \complexityplot[color=color4,mark=x]{cap_None_False}{data/uniform/0.0025_160000.csv}
    \complexityplot[color=color1,mark=o]{cap_None_True}{data/uniform/0.0025_160000.csv}
    
    \nextgroupplot[
        legend style={
            at={(1.0,1.02)},anchor=south east,
        },
        legend columns=5,
        font=\footnotesize,
        yticklabel pos=right,
    ]

    \complexityplot[color=color2,mark=square, solid]{qsvt_None}{data/uniform/0.005_160000.csv}
    \complexityplot[color=color3,mark=triangle, solid]{q_cheb_None}{data/uniform/0.005_160000.csv}
    \complexityplot[color=color3,mark=pentagon, dashed]{chebopt_None}{data/uniform/0.005_160000.csv}
    \complexityplot[color=color4,mark=x]{cap_None_False}{data/uniform/0.005_160000.csv}
    \complexityplot[color=color1,mark=o]{cap_None_True}{data/uniform/0.005_160000.csv}

    \addlegendentry{$\qsvt$}
    \addlegendentry{$\qcheb$}
    \addlegendentry{$\chebopt$}
    \addlegendentry{$\cups$}
    \addlegendentry{$\caps$}
    \nextgroupplot[
        ylabel={$\err(A, b)$ clustered $\lambda$},
        font=\footnotesize,
    ]

    \complexityplot[color=color2,mark=square, solid]{qsvt_None}{data/clustered/0_160000.csv}
    \complexityplot[color=color3,mark=triangle, solid]{q_cheb_None}{data/clustered/0_160000.csv}
    \complexityplot[color=color3,mark=pentagon, dashed]{chebopt_None}{data/clustered/0_160000.csv}
    \complexityplot[color=color4,mark=x]{cap_None_False}{data/clustered/0_160000.csv}
    \complexityplot[color=color1,mark=o]{cap_None_True}{data/clustered/0_160000.csv}
    
    \nextgroupplot[
        font=\footnotesize,
        xlabel={complexity},
        yticklabels=\empty,
    ]
    
    \complexityplot[color=color2,mark=square, solid]{qsvt_None}{data/clustered/0.0025_160000.csv}
    \complexityplot[color=color3,mark=triangle, solid]{q_cheb_None}{data/clustered/0.0025_160000.csv}
    \complexityplot[color=color3,mark=pentagon, dashed]{chebopt_None}{data/clustered/0.0025_160000.csv}
    \complexityplot[color=color4,mark=x]{cap_None_False}{data/clustered/0.0025_160000.csv}
    \complexityplot[color=color1,mark=o]{cap_None_True}{data/clustered/0.0025_160000.csv}
    
    \nextgroupplot[
        legend style={
            at={(1.0,1.02)},anchor=south east,
        },
        legend columns=4,
        font=\footnotesize,
        yticklabel pos=right,
    ]

    \complexityplot[color=color2,mark=square, solid]{qsvt_None}{data/clustered/0.005_160000.csv}
    \complexityplot[color=color3,mark=triangle, solid]{q_cheb_None}{data/clustered/0.005_160000.csv}
    \complexityplot[color=color3,mark=pentagon, dashed]{chebopt_None}{data/clustered/0.005_160000.csv}
    \complexityplot[color=color4,mark=x]{cap_None_False}{data/clustered/0.005_160000.csv}
    \complexityplot[color=color1,mark=o]{cap_None_True}{data/clustered/0.005_160000.csv}

    \end{groupplot}
    \draw (group c1r1.north west) node[draw, anchor=north west, outer sep=0.2cm] {$\xi = 0$};
    \draw (group c2r1.north west) node[draw, anchor=north west, outer sep=0.2cm] {$\xi = 0.0025$};
    \draw (group c3r1.north west) node[draw, anchor=north west, outer sep=0.2cm] {$\xi = 0.005$};
\end{tikzpicture}
\caption{
Relationship of median error and complexity for the $\qsvt$, $\qcheb$, $\chebopt$, $\cups$, and $\caps$ solvers, with noise parameter~$\xi = 0$ (left), $\xi = 0.0025$ (middle), and $\xi = 0.005$ (right), as well as $\samples = 16 \cdot 10^4$ samples per measurement. See \cref{def:complexity} for the exact definition of complexity and error. Methods are tested on randomly generated equations with condition number $\tilde{\kappa} = 3$.
\label{fig:error-complexity}
}
\end{figure}
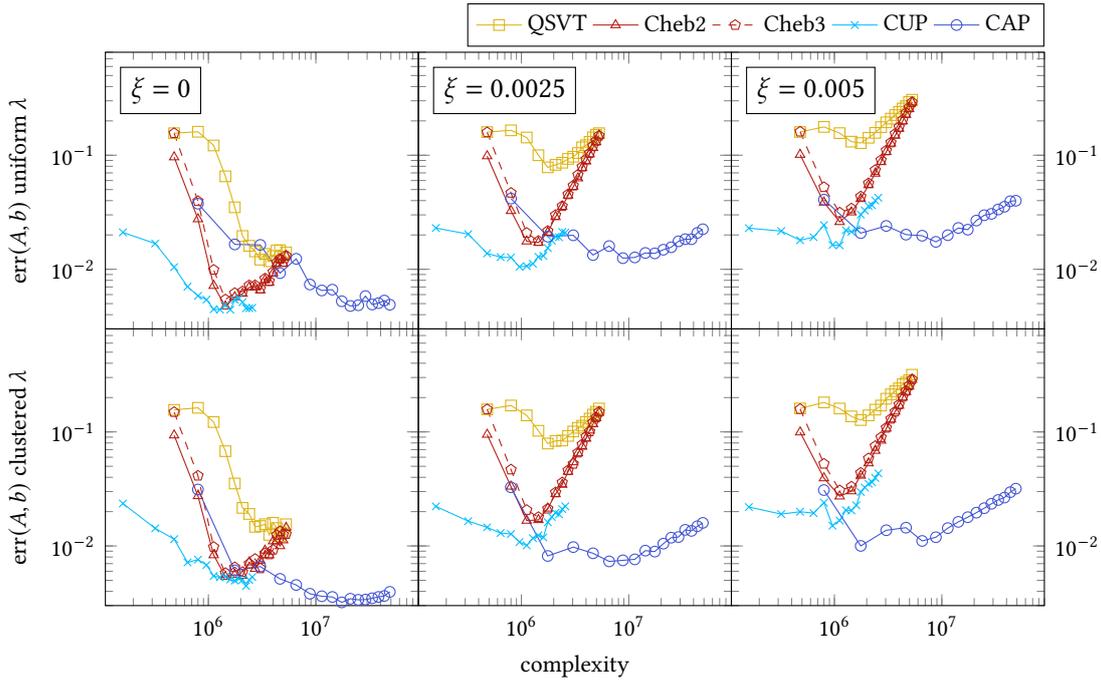

Note that our analysis focuses on the smallest error achievable in the presence of noise. We adopt this perspective because, on current hardware, runtime is typically not the main limitation. Rather, the central challenge is to obtain results that approximate the true solution with meaningful accuracy. For completeness, the tradeoff between error and complexity is shown in \cref{fig:error-complexity}.

\Cref{fig:results} contains convergence plots of the $\qsvt$, $\qcheb$, $\chebopt$, $\cups$, and $\caps$ solvers with different noise rates and without either transformation~\eqref{eq:square-transform} or \eqref{eq:square-outer-transform}.
Note that all five methods do not converge as the polynomial degree increases, even in the noiseless case. This is because the measurement error dominates, which could be decreased by increasing $\samples$. The asymptotic increase of error in the noiseless case is linked to an increase in the sup-norm of the solver polynomials, which could be addressed by multiplying the solver polynomial by a smoothed rectangle function that suppresses it outside the spectral interval, as suggested in \cite{GSLW19}. However, this rectangle function causes a significant overhead, adding approximately $50$ to the total polynomial degree in the current setting. Thus, it would not help in the more realistic noisy simulations, which show an even earlier increase in error. As seen in \cref{tab:all-results}, this noise-induced error is almost independent of the number of samples. With $n$ applications of the block encoding and expected $\xi$ flips per application, the expected total flip count is $\xi\,n$. This scaling is visible in the asymptotic behavior of $\qsvt$ and $\qcheb$ in \cref{fig:results}.

Comparing the five methods, $\cups$ and $\caps$ consistently achieve the lowest error. $\cups$ is slightly better on uniform spectra, while $\caps$ significantly outperforms $\cups$ on clustered spectra. However, $\caps$ is more sensitive than $\cups$ to the over-selection of $n$. Above the optimum, its adaptivity amplifies measurement and noise-induced errors, producing a steeper error rise. At the $95^{\mathrm{th}}$
percentile, the two methods perform similarly, indicating that $\caps$'s advantage lies in typical-case rather than worst-case performance.

\subsection{Hardware results}
To demonstrate that the method works on real hardware, we applied $\caps_\mathrm{sq}$ to the benchmark problem from~\cite[Section~7.2]{DP25}, a linear system arising from the preconditioned discretization of a prototypical elliptic PDE. The right-hand side is strongly aligned with the second-smallest eigenvector of the preconditioned system, producing a spectral weight distribution that is sharply concentrated at one eigenvalue. This is exactly the regime where the adaptivity of $\caps$ is expected to provide the largest gain over non-adaptive methods. Indeed, if $\tilde{\kappa} = \kappa$ is used and $\mu_1$ is measured exactly -- note that $\mu_0 = 1$ automatically follows from the normalization of $b$ -- the degree-1 $\caps$ polynomial becomes approximately tangent to $1/y$ at the dominant eigenvalue, yielding a~near-optimal approximation using only a single application of the block encoding. Moreover, the quantity of interest $b^\top x = b^\top A^{-1}b$ is a moment of the form in~\cref{ssec:bending}, allowing the circuit depth to be halved.

The experiment was run on \texttt{ibm\_aachen} with $n = 1$ and $N = 10^4$ measurement shots per moment. Four independent runs produced relative errors of $24.4\%$, $0.6\%$, $13.6\%$, and $16.1\%$. While individual runs vary substantially, averaging the four solution estimates before computing the error yields only $1.5\%$, compared with the relative error of $18.3\%$ reported in~\cite[Section~7.2]{DP25} for the QSVT solver on the same problem. This indicates that $\caps$ is approximately unbiased under hardware noise, while individual runs have nontrivial variance. This behavior is consistent with the simulated results in \cref{fig:results}.

\section{Conclusion}
Our analysis shows that, although much of the recent progress in quantum linear system solvers is driven by asymptotic complexity, and although Chebyshev-based polynomial constructions demonstrate that better polynomial choices within QSVT meaningfully improve practical performance, there remains considerable room for improving their practical performance by incorporating well-known ideas from classical linear solvers. Within the framework of constrained optimal polynomials developed here, we introduced two classes of solvers: the Constrained Uniform Polynomial (CUP) solver and the Constrained Adaptive Polynomial (CAP) solver. They reduce the noise-induced errors substantially compared with the existing solver. In applications, such an improvement may determine whether quantum simulations yield results that are practically meaningful or not.

In the noise-limited regime, the CUP and CAP solvers outperform the established polynomial constructions while retaining structural advantages. CAP provides a posteriori error estimates and removes the need to accurately prescribe spectral bounds as separate tuning parameters. A~further practical observation concerns the relationship between the two. CAP attains the lower median error for favorable problems but amplifies errors more strongly when the noise passes a~certain threshold, a side effect of its spectral adaptivity. CUP, being problem-independent beyond spectral bounds, is less sensitive to such noise and therefore offers a robust default, while CAP is preferable when the spectrum is known to exhibit structure that adaptivity can exploit. In all cases, the error as a function of polynomial degree is U-shaped; i.e., each solver has an optimal degree beyond which noise dominates, and constrained polynomials lower both the location and the depth of this minimum.

The transition from $\qsvt$ to CUP identifies a simple improvement that benefits not only the semi-iterative solvers studied here but potentially also Variable-Time Amplitude Amplification-based solvers~\cite{Amb10,LS26}, which use semi-iterative solvers as building blocks. An interesting direction for future work is to investigate whether the asymptotic advantages of VTAA or eigenfiltering-based methods can be combined with the adaptive solver framework developed here. More broadly, constrained optimal polynomial design is not specific to matrix inversion. The same perspective applies to other spectral transformations arising as subroutines in quantum algorithms, and extending it to those settings is a natural next step.

\bibliographystyle{quantum}
\bibliography{references}

\vfill

\appendix
\section{Additional data}

This appendix contains a full overview of the best achieved median error using each combination of method and parameters.
The subscripts $\ast_\mathrm{sq}$ and $\ast_\mathrm{sq}'$ indicate \eqref{eq:square-transform} and \eqref{eq:square-outer-transform} respectively. The reported quantity is the median of $\err(A, b)$ over all test equations. For each combination of solver, noise, and samples, the best number of steps~$n$ is selected, indicated in parenthesis after the result. The best solver in each column is highlighted.

\begin{table}[p]
    \centering
    \scalebox{0.87}{
    \begin{tabular}{llllllllll}
        \toprule
        &\multicolumn{3}{c}{$\xi = 0$} & \multicolumn{3}{c}{$\xi = 0.0025$} & \multicolumn{3}{c}{$\xi = 0.005$} \\ \cmidrule(l){2-4}\cmidrule(l){5-7}\cmidrule(l){8-10}
        solver & $\samples = 10^4$ & $4 \cdot 10^4$ & $16 \cdot 10^4$ & $10^4$ & $4 \cdot 10^4$ & $16 \cdot 10^4$ & $10^4$ & $4 \cdot 10^4$ & $16 \cdot 10^4$ \\
        \midrule
$\qsvt$ & 0.038(9) & 0.021(9) & 0.012(11) & 0.065(7) & 0.070(5) & 0.078(5) & 0.108(5) & 0.118(5) & 0.128(5) \\
$\cheb$ & 0.091(1) & 0.048(1) & 0.022(1) & 0.094(1) & 0.045(1) & 0.025(1) & 0.079(1) & 0.040(1) & 0.022(1) \\
$\qcheb$ & 0.016(3) & 0.009(3) & 0.005(4) & 0.022(3) & 0.018(3) & 0.017(4) & 0.026(3) & 0.025(3) & 0.026(3) \\
$\chebopt$ & 0.018(3) & 0.010(4) & 0.005(4) & 0.023(4) & 0.017(4) & 0.018(4) & 0.031(3) & 0.033(3) & 0.031(3) \\
$\cups$ & \textbf{0.014(7)} & \textbf{0.007(7)} & 0.004(10) & \textbf{0.018(7)} & \textbf{0.013(7)} & \textbf{0.011(6)} & \textbf{0.022(7)} & \textbf{0.016(7)} & \textbf{0.016(7)} \\
$\caps$ & 0.018(9) & 0.009(9) & 0.005(10) & 0.022(5) & 0.017(10) & 0.013(6) & 0.026(5) & 0.022(2) & 0.017(6) \\
$\qsvt_\mathrm{sq}$ & 0.038(14) & 0.021(14) & 0.012(14) & 0.105(6) & 0.117(5) & 0.128(5) & 0.151(1) & 0.158(1) & 0.164(1) \\
$\cheb_\mathrm{sq}$ & 0.029(1) & 0.016(1) & 0.011(1) & 0.032(1) & 0.020(1) & 0.016(1) & 0.027(1) & 0.018(1) & 0.018(1) \\
$\qcheb_\mathrm{sq}$ & 0.018(3) & 0.009(4) & 0.005(4) & 0.027(3) & 0.028(3) & 0.027(3) & 0.046(3) & 0.047(3) & 0.047(3) \\
$\chebopt_\mathrm{sq}$ & 0.020(4) & 0.010(4) & 0.005(5) & 0.033(4) & 0.032(3) & 0.032(3) & 0.055(3) & 0.052(3) & 0.052(3) \\
$\cups_\mathrm{sq}$ & 0.018(13) & 0.009(10) & 0.004(10) & 0.028(8) & 0.018(1) & 0.015(1) & 0.025(1) & 0.017(1) & 0.016(1) \\
$\caps_\mathrm{sq}$ & 0.020(16) & 0.009(16) & 0.004(16) & 0.025(6) & 0.014(1) & 0.012(1) & 0.029(1) & 0.020(1) & 0.018(1) \\
$\qsvt_\mathrm{sq}'$ & 0.071(9) & 0.043(9) & 0.027(9) & 0.144(3) & 0.179(3) & 0.219(3) & 0.252(3) & 0.311(2) & 0.366(2) \\
$\qcheb_\mathrm{sq}'$ & 0.019(2) & 0.010(2) & 0.005(2) & 0.024(2) & 0.020(2) & 0.019(2) & 0.037(2) & 0.036(2) & 0.036(2) \\
$\chebopt_\mathrm{sq}'$ & 0.020(2) & 0.011(2) & 0.006(2) & 0.026(2) & 0.021(2) & 0.021(2) & 0.044(2) & 0.041(2) & 0.042(2) \\
$\cups_\mathrm{sq}'$ & 0.016(12) & 0.008(11) & 0.004(11) & 0.022(5) & 0.018(5) & 0.018(6) & 0.034(6) & 0.033(6) & 0.033(6) \\
$\caps_\mathrm{sq}'$ & 0.015(5) & 0.008(5) & \textbf{0.004(15)} & 0.018(6) & 0.017(5) & 0.018(6) & 0.033(6) & 0.033(5) & 0.033(6) \\
        \bottomrule
        \toprule
        &\multicolumn{3}{c}{$\xi = 0.01$} & \multicolumn{3}{c}{$\xi = 0.02$} & \multicolumn{3}{c}{$\xi = 0.04$} \\ \cmidrule(l){2-4}\cmidrule(l){5-7}\cmidrule(l){8-10}
        solver & $\samples = 10^4$ & $4 \cdot 10^4$ & $16 \cdot 10^4$ & $10^4$ & $4 \cdot 10^4$ & $16 \cdot 10^4$ & $10^4$ & $4 \cdot 10^4$ & $16 \cdot 10^4$ \\
        \midrule
$\qsvt$ & 0.150(1) & 0.159(1) & 0.161(1) & 0.159(1) & 0.166(1) & 0.171(1) & 0.175(1) & 0.181(1) & 0.186(1) \\
$\cheb$ & 0.084(1) & 0.044(1) & 0.022(1) & 0.093(1) & 0.049(1) & \textbf{0.024(1)} & 0.087(1) & 0.044(1) & \textbf{0.033(1)} \\
$\qcheb$ & 0.049(2) & 0.048(3) & 0.047(3) & 0.074(2) & 0.074(2) & 0.073(2) & 0.119(2) & 0.119(2) & 0.118(2) \\
$\chebopt$ & 0.053(3) & 0.052(3) & 0.052(3) & 0.088(2) & 0.089(2) & 0.089(2) & 0.141(2) & 0.141(2) & 0.140(2) \\
$\cups$ & \textbf{0.028(2)} & 0.027(2) & \textbf{0.021(1)} & \textbf{0.039(2)} & 0.038(1) & 0.024(1) & 0.061(2) & \textbf{0.042(1)} & 0.034(1) \\
$\caps$ & 0.032(2) & 0.028(2) & 0.027(2) & 0.041(2) & \textbf{0.036(2)} & 0.036(2) & \textbf{0.059(2)} & 0.059(2) & 0.050(1) \\
$\qsvt_\mathrm{sq}$ & 0.160(1) & 0.165(1) & 0.170(1) & 0.177(1) & 0.183(1) & 0.186(1) & 0.212(1) & 0.216(1) & 0.221(1) \\
$\cheb_\mathrm{sq}$ & 0.031(1) & 0.028(1) & 0.027(1) & 0.048(1) & 0.046(1) & 0.047(1) & 0.091(1) & 0.088(1) & 0.087(1) \\
$\qcheb_\mathrm{sq}$ & 0.071(2) & 0.072(2) & 0.073(2) & 0.115(2) & 0.119(2) & 0.118(2) & 0.173(1) & 0.173(1) & 0.172(1) \\
$\chebopt_\mathrm{sq}$ & 0.090(2) & 0.090(2) & 0.089(2) & 0.143(2) & 0.141(2) & 0.141(2) & 0.226(1) & 0.227(1) & 0.224(1) \\
$\cups_\mathrm{sq}$ & 0.030(1) & \textbf{0.025(1)} & 0.025(1) & 0.045(1) & 0.045(1) & 0.044(1) & 0.086(1) & 0.084(1) & 0.082(1) \\
$\caps_\mathrm{sq}$ & 0.032(1) & 0.028(1) & 0.026(1) & 0.048(1) & 0.044(1) & 0.044(1) & 0.079(1) & 0.081(1) & 0.080(1) \\
$\qsvt_\mathrm{sq}'$ & 0.373(2) & 0.449(2) & 0.530(2) & 0.591(2) & 0.718(2) & 0.838(2) & 0.715(1) & 0.849(1) & 0.966(1) \\
$\qcheb_\mathrm{sq}'$ & 0.062(1) & 0.060(1) & 0.060(1) & 0.083(1) & 0.083(1) & 0.083(1) & 0.128(1) & 0.129(1) & 0.128(1) \\
$\chebopt_\mathrm{sq}'$ & 0.076(1) & 0.076(2) & 0.076(2) & 0.101(1) & 0.101(1) & 0.101(1) & 0.152(1) & 0.153(1) & 0.151(1) \\
$\cups_\mathrm{sq}'$ & 0.055(4) & 0.054(3) & 0.053(3) & 0.079(4) & 0.077(4) & 0.077(3) & 0.125(3) & 0.125(3) & 0.124(3) \\
$\caps_\mathrm{sq}'$ & 0.048(4) & 0.049(3) & 0.050(3) & 0.074(3) & 0.075(3) & 0.074(4) & 0.120(3) & 0.122(4) & 0.123(4) \\
    \bottomrule
\end{tabular}
    }
    \caption{
    Best results for each solver configuration for uniform eigenvalues. 
    }
    \label{tab:all-results}
\end{table}

\begin{table}[p]
    \centering
    \scalebox{0.87}{
    \begin{tabular}{llllllllll}
        \toprule
        &\multicolumn{3}{c}{$\xi = 0$} & \multicolumn{3}{c}{$\xi = 0.0025$} & \multicolumn{3}{c}{$\xi = 0.005$} \\ \cmidrule(l){2-4}\cmidrule(l){5-7}\cmidrule(l){8-10}
        solver & $\samples = 10^4$ & $4 \cdot 10^4$ & $16 \cdot 10^4$ & $10^4$ & $4 \cdot 10^4$ & $16 \cdot 10^4$ & $10^4$ & $4 \cdot 10^4$ & $16 \cdot 10^4$ \\
        \midrule
$\qsvt$ & 0.040(14) & 0.022(14) & 0.012(14) & 0.066(5) & 0.069(5) & 0.079(5) & 0.103(5) & 0.116(6) & 0.127(5) \\
$\cheb$ & 0.089(1) & 0.043(1) & 0.024(1) & 0.087(1) & 0.043(1) & 0.023(1) & 0.084(1) & 0.040(1) & 0.022(1) \\
$\qcheb$ & 0.018(3) & 0.010(4) & 0.005(4) & 0.022(3) & 0.017(3) & 0.017(3) & 0.029(3) & 0.028(3) & 0.027(3) \\
$\chebopt$ & 0.019(3) & 0.010(4) & 0.006(4) & 0.024(4) & 0.018(4) & 0.018(4) & 0.034(4) & 0.031(3) & 0.031(3) \\
$\cups$ & 0.016(14) & 0.009(14) & 0.004(14) & 0.020(5) & 0.012(7) & 0.010(7) & 0.023(5) & 0.016(6) & 0.015(6) \\
$\caps$ & 0.014(5) & \textbf{0.006(9)} & 0.003(9) & \textbf{0.016(5)} & \textbf{0.010(5)} & \textbf{0.007(5)} & \textbf{0.017(2)} & \textbf{0.013(2)} & \textbf{0.010(2)} \\
$\qsvt_\mathrm{sq}$ & 0.040(15) & 0.023(15) & 0.013(15) & 0.104(5) & 0.111(5) & 0.125(5) & 0.150(3) & 0.157(1) & 0.165(1) \\
$\cheb_\mathrm{sq}$ & 0.032(1) & 0.017(1) & 0.012(1) & 0.031(1) & 0.020(1) & 0.015(1) & 0.030(1) & 0.020(1) & 0.020(1) \\
$\qcheb_\mathrm{sq}$ & 0.018(3) & 0.009(4) & 0.005(4) & 0.032(3) & 0.029(3) & 0.027(3) & 0.044(3) & 0.046(3) & 0.047(3) \\
$\chebopt_\mathrm{sq}$ & 0.022(3) & 0.011(6) & 0.006(5) & 0.030(4) & 0.031(3) & 0.032(3) & 0.056(3) & 0.053(3) & 0.052(3) \\
$\cups_\mathrm{sq}$ & 0.019(15) & 0.009(15) & 0.004(15) & 0.030(1) & 0.019(1) & 0.014(1) & 0.026(1) & 0.020(1) & 0.018(1) \\
$\caps_\mathrm{sq}$ & 0.014(13) & 0.007(9) & \textbf{0.003(9)} & 0.024(6) & 0.016(1) & 0.011(1) & 0.028(1) & 0.016(1) & 0.014(1) \\
$\qsvt_\mathrm{sq}'$ & 0.078(7) & 0.050(7) & 0.031(7) & 0.159(4) & 0.181(3) & 0.224(3) & 0.242(3) & 0.311(3) & 0.381(2) \\
$\qcheb_\mathrm{sq}'$ & 0.020(2) & 0.010(2) & 0.006(2) & 0.025(2) & 0.020(2) & 0.020(2) & 0.035(2) & 0.036(2) & 0.035(2) \\
$\chebopt_\mathrm{sq}'$ & 0.021(2) & 0.011(2) & 0.007(2) & 0.029(2) & 0.022(2) & 0.021(2) & 0.044(2) & 0.042(2) & 0.042(2) \\
$\cups_\mathrm{sq}'$ & 0.017(5) & 0.009(6) & 0.005(11) & 0.022(5) & 0.017(6) & 0.018(5) & 0.034(6) & 0.032(6) & 0.032(6) \\
$\caps_\mathrm{sq}'$ & \textbf{0.013(10)} & 0.006(14) & 0.003(14) & 0.018(4) & 0.015(6) & 0.013(4) & 0.024(3) & 0.020(4) & 0.020(3) \\
        \bottomrule
        \toprule
        &\multicolumn{3}{c}{$\xi = 0.01$} & \multicolumn{3}{c}{$\xi = 0.02$} & \multicolumn{3}{c}{$\xi = 0.04$} \\ \cmidrule(l){2-4}\cmidrule(l){5-7}\cmidrule(l){8-10}
        solver & $\samples = 10^4$ & $4 \cdot 10^4$ & $16 \cdot 10^4$ & $10^4$ & $4 \cdot 10^4$ & $16 \cdot 10^4$ & $10^4$ & $4 \cdot 10^4$ & $16 \cdot 10^4$ \\
        \midrule
$\qsvt$ & 0.151(3) & 0.158(1) & 0.164(1) & 0.160(1) & 0.174(1) & 0.176(1) & 0.174(1) & 0.184(1) & 0.193(1) \\
$\cheb$ & 0.078(1) & 0.041(1) & 0.023(1) & 0.085(1) & 0.046(1) & \textbf{0.025(1)} & 0.095(1) & 0.056(1) & 0.039(1) \\
$\qcheb$ & 0.043(3) & 0.044(3) & 0.045(3) & 0.073(2) & 0.075(2) & 0.074(2) & 0.120(2) & 0.121(2) & 0.121(2) \\
$\chebopt$ & 0.051(3) & 0.052(3) & 0.050(3) & 0.090(2) & 0.088(2) & 0.089(2) & 0.137(2) & 0.142(2) & 0.141(2) \\
$\cups$ & 0.029(5) & 0.027(2) & 0.023(1) & 0.039(2) & 0.036(2) & 0.026(1) & 0.061(2) & 0.047(1) & \textbf{0.037(1)} \\
$\caps$ & \textbf{0.022(2)} & \textbf{0.017(2)} & \textbf{0.017(2)} & \textbf{0.031(2)} & \textbf{0.029(2)} & 0.028(2) & \textbf{0.051(2)} & \textbf{0.047(2)} & 0.041(1) \\
$\qsvt_\mathrm{sq}$ & 0.160(1) & 0.167(1) & 0.171(1) & 0.179(1) & 0.189(1) & 0.189(1) & 0.215(1) & 0.222(1) & 0.224(1) \\
$\cheb_\mathrm{sq}$ & 0.036(1) & 0.028(1) & 0.028(1) & 0.049(1) & 0.047(1) & 0.047(1) & 0.085(1) & 0.083(1) & 0.084(1) \\
$\qcheb_\mathrm{sq}$ & 0.073(2) & 0.073(2) & 0.074(2) & 0.121(2) & 0.120(2) & 0.120(2) & 0.174(1) & 0.172(1) & 0.172(1) \\
$\chebopt_\mathrm{sq}$ & 0.091(2) & 0.092(2) & 0.092(3) & 0.138(2) & 0.137(2) & 0.141(2) & 0.223(1) & 0.223(1) & 0.224(1) \\
$\cups_\mathrm{sq}$ & 0.034(1) & 0.026(1) & 0.025(1) & 0.045(1) & 0.045(1) & 0.044(1) & 0.081(1) & 0.080(1) & 0.079(1) \\
$\caps_\mathrm{sq}$ & 0.033(1) & 0.026(1) & 0.023(1) & 0.047(1) & 0.043(1) & 0.042(1) & 0.083(1) & 0.080(1) & 0.080(1) \\
$\qsvt_\mathrm{sq}'$ & 0.383(2) & 0.458(2) & 0.538(2) & 0.606(2) & 0.713(2) & 0.833(2) & 0.763(1) & 0.901(1) & 1.015(1) \\
$\qcheb_\mathrm{sq}'$ & 0.060(1) & 0.060(1) & 0.060(1) & 0.084(1) & 0.084(1) & 0.083(1) & 0.133(1) & 0.132(1) & 0.132(1) \\
$\chebopt_\mathrm{sq}'$ & 0.074(1) & 0.076(1) & 0.076(2) & 0.101(1) & 0.100(1) & 0.100(1) & 0.151(1) & 0.149(1) & 0.148(1) \\
$\cups_\mathrm{sq}'$ & 0.053(3) & 0.051(3) & 0.051(4) & 0.077(3) & 0.075(3) & 0.074(3) & 0.127(4) & 0.126(3) & 0.125(4) \\
$\caps_\mathrm{sq}'$ & 0.034(3) & 0.034(4) & 0.032(4) & 0.057(4) & 0.056(4) & 0.058(3) & 0.108(4) & 0.108(4) & 0.108(4) \\
    \bottomrule
\end{tabular}
    }
    \caption{
    Best results for each solver configuration for clustered eigenvalues.
    }
    \label{tab:all-results-clustered}
\end{table}

\end{document}

%% file: preamble.tex
\usepackage[utf8]{inputenc}
\usepackage[english]{babel}
\usepackage[T1]{fontenc}
\usepackage{microtype}
\usepackage[sb,rm,tt=false]{libertine}
\usepackage[lf,medium]{FiraSans}
\usepackage[T1]{fontenc}
\usepackage{textcomp}
\usepackage[varqu,varl]{zi4}
\usepackage{mathtools}
\usepackage[amsthm,libertine,vvarbb]{newtxmath}
\usepackage{thmtools}
\usepackage{color}
\usepackage{pgfplots}
\pgfplotsset{compat=1.18}
\usepgfplotslibrary{groupplots,fillbetween}
\usepackage{tikz}
\usetikzlibrary{quantikz2,external}

\usepackage[numbers]{natbib}
\usepackage[hyphens]{url}
\usepackage[hidelinks]{hyperref}
\usepackage{titlesec}
\usepackage[ruled]{algorithm2e}
\usepackage[font=footnotesize]{subcaption}
\usepackage{enumitem}
\usepackage{mathdots}
\usepackage[nameinlink,noabbrev,capitalise]{cleveref}
\usepackage{booktabs}

\DisableLigatures[q]{shape = sc}

\newtheorem{theorem}{Theorem}

\newtheorem{lemma}[theorem]{Lemma}
\theoremstyle{definition}
\newtheorem{definition}[theorem]{Definition}

\theoremstyle{remark}

\Crefname{assumption}{Assumption}{Assumptions}

\setlist[enumerate]{label=\arabic*.}

\definecolor{unia-mntf}{RGB}{0, 101, 97}
\definecolor{unia-purple}{RGB}{173, 0, 124}

\hypersetup{
    colorlinks=true,
    linkcolor=unia-purple,
    urlcolor=unia-purple,
    citecolor=unia-purple
}

\DeclareFontFamily{U}{matha}{\hyphenchar\font45}
\DeclareFontShape{U}{matha}{m}{n}{
	<-6> matha5 <6-7> matha6 <7-8> matha7
	<8-9> matha8 <9-10> matha9
	<10-12> matha10 <12-> matha12
}{}
\DeclareSymbolFont{matha}{U}{matha}{m}{n}

\DeclareFontFamily{U}{mathx}{\hyphenchar\font45}
\DeclareFontShape{U}{mathx}{m}{n}{
	<-6> mathx5 <6-7> mathx6 <7-8> mathx7
	<8-9> mathx8 <9-10> mathx9
	<10-12> mathx10 <12-> mathx12
}{}
\DeclareSymbolFont{mathx}{U}{mathx}{m}{n}

\DeclareMathDelimiter{\vvvert} {0}{matha}{"7E}{mathx}{"17}

\allowdisplaybreaks

\makeatletter
\def\th@plain{%
  \thm@headfont{\scshape}
  \thm@notefont{\normalfont}% same as heading font
  \itshape % body font
}
\def\th@definition{%
  \thm@headfont{\scshape}
  \thm@notefont{\normalfont}% same as heading font
  \normalfont % body font
}
\makeatother

\SetKwSty{textsc}

\makeatletter
\def\@setauthors{%
  \begingroup
  \def\thanks{\protect\thanks@warning}%
  \trivlist
  \raggedright\footnotesize \@topsep30\p@\relax
  \advance\@topsep by -\baselineskip
  \item\relax
  \author@andify\authors
  \def\\{\protect\linebreak}%
  \rmfamily\large\authors%
  \ifx\@empty\contribs
  \else
    ,\penalty-3 \space \@setcontribs
    \@closetoccontribs
  \fi
  \endtrivlist
  \endgroup
}
\def\@settitle{%
    \noindent
    \baselineskip14\p@\relax
    \LARGE
    \firabook\sffamily\raggedright
    \@title
}
\def\@setaddresses{\par
  \nobreak \begingroup
  \def\author##1{\nobreak\addvspace\smallskipamount}%
  \def\\{\unskip, \ignorespaces}%
  \interlinepenalty\@M%
  \def\address##1##2{\begingroup%
    \par\addvspace\smallskipamount\indent%
    \@ifnotempty{##1}{(\ignorespaces##1\unskip) }%
    {\ignorespaces##2}\par\endgroup}%
  \def\curraddr##1##2{\begingroup%
    \@ifnotempty{##2}{\nobreak\indent\curraddrname
      \@ifnotempty{##1}{, \ignorespaces##1\unskip}\/:\space
      ##2\par}\endgroup}%
  \def\email##1##2{\begingroup
    \@ifnotempty{##2}{\nobreak\indent\emailaddrname
      \@ifnotempty{##1}{, \ignorespaces##1\unskip}\/:\space
      \ttfamily##2\par}\endgroup}%
  \def\urladdr##1##2{\begingroup
    \def~{\char`\~}%
    \@ifnotempty{##2}{\nobreak\indent\urladdrname
      \@ifnotempty{##1}{, \ignorespaces##1\unskip}\/:\space
      \ttfamily##2\par}\endgroup}%
\begin{flushleft}
    \addressfont\addresses
  \end{flushleft}
  \endgroup
}
\newcommand{\addressfont}{\normalfont\small}
\patchcmd{\abstract}{3pc}{0pt}{}{} % remove indentation
\patchcmd{\abstract}{\parsep}{\setlength{\itemsep}{0.8em}\parsep}{}{} % remove indentation
\renewcommand{\@mkboth}[2]{} % no automatic running head
\patchcmd{\@maketitle}
  {\ifx\@empty\@dedicatory}
  {\@setaddresses\global\let\@setaddresses\relax\ifx\@empty\@dedicatory}
  {}{}
\makeatother

\titleformat{\section}
  {\sffamily\firabook\large} % Commands for section titles
  {\thesection}
  {1em}
  {}
\titleformat{\subsection}
  {\sffamily\firabook\normalsize} % Commands for section titles
  {\thesubsection}
  {1em}
  {}
\titlespacing*{\section} {0pt}{3.5ex plus 1ex minus .2ex}{2.3ex plus .2ex}
\titlespacing*{\subsection} {0pt}{3.25ex plus 1ex minus .2ex}{1.5ex plus .2ex}

\titleformat{\paragraph}[runin]{\normalfont\normalsize\itshape}{\theparagraph}{1em}{}[.]
\titlespacing*{\paragraph}{0pt}{3.25ex plus 1ex minus .2ex}{\the\fontdimen2\font}

%% file: references.bib
@inproceedings{AJL06,
  title = {A Polynomial Quantum Algorithm for Approximating the {{Jones}} Polynomial},
  booktitle = {Proc. {{Thirty-Eighth Annu}}. {{ACM Symp}}. {{Theory Comput}}.},
  year = 2006,
  month = may,
  eprint = {quant-ph/0511096},
  pages = {427--436},
  publisher = {ACM},
  address = {Seattle WA USA},
  doi = {10.1145/1132516.1132579},
  urldate = {2026-03-23},
  archiveprefix = {arXiv},
  isbn = {978-1-59593-134-4},
  langid = {english},
  author = {Aharonov, D. and Jones, V. and Landau, Z.}
}

@article{AL22,
  title = {Quantum Linear System Solver Based on Time-Optimal Adiabatic Quantum Computing and Quantum Approximate Optimization Algorithm},
  year = 2022,
  month = jun,
  journal = {ACM Transactions on Quantum Computing},
  volume = {3},
  number = {2},
  eprint = {1909.05500},
  primaryclass = {quant-ph},
  pages = {1--28},
  issn = {2643-6809, 2643-6817},
  doi = {10.1145/3498331},
  urldate = {2023-10-12},
  archiveprefix = {arXiv},
  langid = {english},
  author = {An, D. and Lin, L.}
}

@article{Amb10,
  title = {Variable Time Amplitude Amplification and a Faster Quantum Algorithm for Solving Systems of Linear Equations},
  year = 2010,
  month = nov,
  number = {arXiv:1010.4458},
  eprint = {1010.4458},
  primaryclass = {quant-ph},
  publisher = {arXiv},
  urldate = {2025-05-30},
  archiveprefix = {arXiv},
  author = {Ambainis, A.}
}

@article{BDP26,
  title = {Quantum {{Enhanced Numerical Homogenization}}},
  year = 2026,
  eprint = {2603.28521},
  publisher = {arXiv},
  urldate = {2026-04-02},
  archiveprefix = {arXiv},
  author = {Balazi, L. and Deiml, M. and Peterseim, D.}
}

@article{CGJ19,
  title = {The {{Power}} of {{Block-Encoded Matrix Powers}}: {{Improved Regression Techniques}} via {{Faster Hamiltonian Simulation}}},
  shorttitle = {The {{Power}} of {{Block-Encoded Matrix Powers}}},
  year = 2019,
  journal = {LIPIcs Vol. 132 ICALP 2019},
  volume = {132},
  eprint = {1804.01973},
  pages = {33:1-33:14},
  publisher = {Schloss Dagstuhl -- Leibniz-Zentrum f\"ur Informatik},
  issn = {1868-8969},
  doi = {10.4230/LIPICS.ICALP.2019.33},
  urldate = {2025-05-30},
  archiveprefix = {arXiv},
  isbn = {9783959771092},
  langid = {english},
  author = {Chakraborty, S. and Gily{\'e}n, A. and Jeffery, S.},
  editor = {Baier, C. and Chatzigiannakis, I. and Flocchini, P. and Leonardi, S.}
}

@article{CKS17,
  title = {Quantum {{Algorithm}} for {{Systems}} of {{Linear Equations}} with {{Exponentially Improved Dependence}} on {{Precision}}},
  year = 2017,
  month = jan,
  journal = {SIAM J. Comput.},
  volume = {46},
  number = {6},
  eprint = {1511.02306},
  pages = {1920--1950},
  issn = {0097-5397, 1095-7111},
  doi = {10.1137/16M1087072},
  urldate = {2024-02-14},
  archiveprefix = {arXiv},
  langid = {english},
  author = {Childs, A. M. and Kothari, R. and Somma, R. D.}
}

@article{Dal24,
  title = {A Shortcut to an Optimal Quantum Linear System Solver},
  year = 2024,
  month = jun,
  number = {arXiv:2406.12086},
  eprint = {2406.12086},
  primaryclass = {quant-ph},
  publisher = {arXiv},
  urldate = {2024-10-28},
  archiveprefix = {arXiv},
  langid = {english},
  author = {Dalzell, A. M.}
}

@article{DMWL21,
  title = {Efficient Phase-Factor Evaluation in Quantum Signal Processing},
  year = 2021,
  month = apr,
  journal = {Phys. Rev. A},
  volume = {103},
  number = {4},
  eprint = {2002.11649},
  primaryclass = {quant-ph},
  pages = {042419},
  issn = {2469-9926, 2469-9934},
  doi = {10.1103/PhysRevA.103.042419},
  urldate = {2025-12-18},
  archiveprefix = {arXiv},
  author = {Dong, Y. and Meng, X. and Whaley, K. B. and Lin, L.}
}

@article{DP24,
  title = {Nonlinear Quantum Computing by Amplified Encodings},
  year = 2024,
  month = nov,
  number = {arXiv:2411.16435},
  eprint = {2411.16435},
  primaryclass = {quant-ph},
  publisher = {arXiv},
  urldate = {2025-03-10},
  archiveprefix = {arXiv},
  author = {Deiml, M. and Peterseim, D.}
}

@article{DP25,
  title = {Quantum {{Realization}} of the {{Finite Element Method}}},
  year = 2025,
  journal = {Math. Comp.},
  eprint = {2403.19512},
  primaryclass = {quant-ph},
  doi = {10.1090/mcom/4124},
  urldate = {2024-04-19},
  archiveprefix = {arXiv},
  langid = {english},
  author = {Deiml, M. and Peterseim, D.}
}

@article{ELN22,
  title = {A {{Theory}} of {{Quantum Subspace Diagonalization}}},
  year = 2022,
  month = sep,
  journal = {SIAM J. Matrix Anal. Appl.},
  volume = {43},
  number = {3},
  eprint = {2110.07492},
  pages = {1263--1290},
  issn = {0895-4798, 1095-7162},
  doi = {10.1137/21M145954X},
  urldate = {2026-04-13},
  archiveprefix = {arXiv},
  langid = {english},
  author = {Epperly, E. N. and Lin, L. and Nakatsukasa, Y.}
}

@article{FS50,
  title = {Numerical {{Determination}} of {{Fundamental Modes}}},
  year = 1950,
  month = dec,
  journal = {J. Appl. Phys.},
  volume = {21},
  number = {12},
  pages = {1326--1332},
  issn = {0021-8979, 1089-7550},
  doi = {10.1063/1.1699598},
  urldate = {2025-12-15},
  langid = {english},
  author = {Flanders, D. A. and Shortley, G.}
}

@article{GKS24,
  title = {An {{Optimal Linear-combination-of-unitaries-based Quantum Linear System Solver}}},
  year = 2024,
  month = jun,
  journal = {ACM Trans. Quantum Comput.},
  volume = {5},
  number = {2},
  pages = {1--23},
  issn = {2643-6809, 2643-6817},
  doi = {10.1145/3649320},
  urldate = {2026-04-24},
  langid = {english},
  author = {Gribling, S. and Kerenidis, I. and Szil{\'a}gyi, D.}
}

@inproceedings{GSLW19,
  title = {Quantum Singular Value Transformation and beyond: Exponential Improvements for Quantum Matrix Arithmetics},
  shorttitle = {Quantum Singular Value Transformation and Beyond},
  booktitle = {Proc. 51st {{Annu}}. {{ACM SIGACT Symp}}. {{Theory Comput}}.},
  year = 2019,
  month = jun,
  eprint = {1806.01838},
  primaryclass = {quant-ph},
  pages = {193--204},
  doi = {10.1145/3313276.3316366},
  urldate = {2023-10-12},
  archiveprefix = {arXiv},
  langid = {english},
  author = {Gily{\'e}n, A. and Su, Y. and Low, G. H. and Wiebe, N.}
}

@article{HHL09,
  title = {Quantum Algorithm for Solving Linear Systems of Equations},
  year = 2009,
  month = oct,
  journal = {Phys. Rev. Lett.},
  volume = {103},
  number = {15},
  eprint = {0811.3171},
  primaryclass = {quant-ph},
  pages = {150502},
  issn = {0031-9007, 1079-7114},
  doi = {10.1103/PhysRevLett.103.150502},
  urldate = {2023-10-24},
  archiveprefix = {arXiv},
  langid = {english},
  author = {Harrow, A. W. and Hassidim, A. and Lloyd, S.}
}

@article{KMM23,
  title = {Exact and Efficient {{Lanczos}} Method on a Quantum Computer},
  year = 2023,
  month = may,
  journal = {Quantum},
  volume = {7},
  eprint = {2208.00567},
  primaryclass = {quant-ph},
  pages = {1018},
  issn = {2521-327X},
  doi = {10.22331/q-2023-05-23-1018},
  urldate = {2025-05-26},
  archiveprefix = {arXiv},
  author = {Kirby, W. and Motta, M. and Mezzacapo, A.}
}

@techreport{Kra88,
  type = {Technical {{Report}}},
  title = {A Software Package for Sequential Quadratic Programming},
  year = 1988,
  month = jul,
  number = {DFVLR-FB 88-28},
  institution = {Institut f\"ur Dynamik der Flugsysteme, Deutsche Forschungs- und Versuchsanstalt f\"ur Luft- und Raumfahrt (DFVLR), Oberpfaffenhofen},
  author = {Kraft, D.}
}

@book{LS12,
  title = {Krylov {{Subspace Methods}}: {{Principles}} and {{Analysis}}},
  shorttitle = {Krylov {{Subspace Methods}}},
  year = 2012,
  month = oct,
  publisher = {Oxford University Press},
  doi = {10.1093/acprof:oso/9780199655410.001.0001},
  urldate = {2026-03-25},
  isbn = {978-0-19-965541-0},
  author = {Liesen, J. and Strakos, Z.}
}

@article{LS26,
  title = {Quantum Linear System Algorithm with Optimal Queries to Initial State Preparation},
  year = 2026,
  month = mar,
  journal = {Quantum},
  volume = {10},
  eprint = {2410.18178},
  primaryclass = {quant-ph},
  pages = {2041},
  issn = {2521-327X},
  doi = {10.22331/q-2026-03-23-2041},
  urldate = {2026-04-02},
  archiveprefix = {arXiv},
  langid = {english},
  author = {Low, G. H. and Su, Y.}
}

@article{LT20,
  title = {Optimal Polynomial Based Quantum Eigenstate Filtering with Application to Solving Quantum Linear Systems},
  year = 2020,
  month = nov,
  journal = {Quantum},
  volume = {4},
  eprint = {1910.14596},
  primaryclass = {quant-ph},
  pages = {361},
  issn = {2521-327X},
  doi = {10.22331/q-2020-11-11-361},
  urldate = {2023-11-09},
  archiveprefix = {arXiv},
  langid = {english},
  author = {Lin, L. and Tong, Y.}
}

@article{LWX25,
  title = {Generalized Quantum Singular Value Transformation with Application in Quantum Bi-Conjugate Gradient Method},
  year = 2025,
  month = aug,
  number = {arXiv:2508.21390},
  eprint = {2508.21390},
  primaryclass = {math},
  publisher = {arXiv},
  urldate = {2025-12-03},
  archiveprefix = {arXiv},
  langid = {english},
  author = {Liu, Y.-Q. and Wang, H. and Xiang, H.}
}

@article{MP84,
  title = {Maximum Entropy in the Problem of Moments},
  year = 1984,
  month = aug,
  journal = {J. Math. Phys.},
  volume = {25},
  number = {8},
  pages = {2404--2417},
  issn = {0022-2488, 1089-7658},
  doi = {10.1063/1.526446},
  urldate = {2026-04-16},
  langid = {english},
  author = {Mead, L. R. and Papanicolaou, N.}
}

@article{MPS+25,
  title = {Quantum {{Linear System Solvers}}: {{A Survey}} of {{Algorithms}} and {{Applications}}},
  shorttitle = {Quantum {{Linear System Solvers}}},
  year = 2025,
  month = jan,
  number = {arXiv:2411.02522},
  eprint = {2411.02522},
  primaryclass = {quant-ph},
  publisher = {arXiv},
  urldate = {2025-06-20},
  archiveprefix = {arXiv},
  author = {Morales, M. E. S. and Pira, L. and Schleich, P. and Koor, K. and Costa, P. C. S. and An, D. and {Aspuru-Guzik}, A. and Lin, L. and Rebentrost, P. and Berry, D. W.}
}

@article{NSYL25,
  title = {Inverse Nonlinear Fast {{Fourier}} Transform on {{SU}}(2) with Applications to Quantum Signal Processing},
  year = 2025,
  month = may,
  number = {arXiv:2505.12615},
  eprint = {2505.12615},
  primaryclass = {quant-ph},
  publisher = {arXiv},
  urldate = {2026-02-12},
  archiveprefix = {arXiv},
  author = {Ni, H. and Sarkar, R. and Ying, L. and Lin, L.}
}

@book{Saa03,
  title = {Iterative {{Methods}} for {{Sparse Linear Systems}}},
  year = 2003,
  month = jan,
  edition = {Second},
  publisher = {{Society for Industrial and Applied Mathematics}},
  doi = {10.1137/1.9780898718003},
  urldate = {2025-12-15},
  isbn = {978-0-89871-534-7 978-0-89871-800-3},
  langid = {english},
  author = {Saad, Y.}
}

@article{SNW+25,
  title = {Matrix Inversion Polynomials for the Quantum Singular Value Transformation},
  year = 2025,
  month = jul,
  number = {arXiv:2507.15537},
  eprint = {2507.15537},
  primaryclass = {quant-ph},
  publisher = {arXiv},
  urldate = {2026-04-24},
  archiveprefix = {arXiv},
  author = {S{\"u}nderhauf, C. and N{\'e}meth, Z. and Walayat, A. and Patterson, A. and Berntson, B. K.}
}

@article{SUR+20,
  title = {Amplitude Estimation without Phase Estimation},
  year = 2020,
  month = feb,
  journal = {Quantum Inf Process},
  volume = {19},
  number = {2},
  eprint = {1904.10246},
  primaryclass = {quant-ph},
  pages = {75},
  issn = {1570-0755, 1573-1332},
  doi = {10.1007/s11128-019-2565-2},
  urldate = {2024-08-27},
  archiveprefix = {arXiv},
  langid = {english},
  author = {Suzuki, Y. and Uno, S. and Raymond, R. and Tanaka, T. and Onodera, T. and Yamamoto, N.}
}

@article{TWYH24,
  title = {Quantum Conjugate Gradient Method Using the Positive-Side Quantum Eigenvalue Transformation},
  year = 2024,
  month = apr,
  number = {arXiv:2404.02713},
  eprint = {2404.02713},
  primaryclass = {quant-ph},
  publisher = {arXiv},
  urldate = {2025-12-03},
  archiveprefix = {arXiv},
  author = {Toyoizumi, K. and Wada, K. and Yamamoto, N. and Hoshino, K.}
}

@article{WZ82,
  title = {A Single Quantum Cannot Be Cloned},
  year = 1982,
  month = oct,
  journal = {Nature},
  volume = {299},
  number = {5886},
  pages = {802--803},
  issn = {0028-0836, 1476-4687},
  doi = {10.1038/299802a0},
  urldate = {2026-01-07},
  langid = {english},
  author = {Wootters, W. K. and Zurek, W. H.}
}

@article{XZZ24,
  title = {Quantum {{Krylov-Subspace Method Based Linear Solver}}},
  year = 2024,
  month = may,
  number = {arXiv:2405.06359},
  eprint = {2405.06359},
  primaryclass = {quant-ph},
  publisher = {arXiv},
  urldate = {2025-12-03},
  archiveprefix = {arXiv},
  langid = {english},
  author = {Xu, R.-B. and Zheng, Z.-J. and Zheng, Z.}
}

@article{ZWXL24,
  title = {Measurement-Efficient Quantum {{Krylov}} Subspace Diagonalisation},
  year = 2024,
  month = aug,
  journal = {Quantum},
  volume = {8},
  eprint = {2301.13353},
  pages = {1438},
  issn = {2521-327X},
  doi = {10.22331/q-2024-08-13-1438},
  urldate = {2026-03-31},
  archiveprefix = {arXiv},
  langid = {english},
  author = {Zhang, Z. and Wang, A. and Xu, X. and Li, Y.}
}

@article{DLT22,
  title = {Ground-{{State Preparation}} and {{Energy Estimation}} on {{Early Fault-Tolerant Quantum Computers}} via {{Quantum Eigenvalue Transformation}} of {{Unitary Matrices}}},
  year = 2022,
  month = oct,
  journal = {PRX Quantum},
  volume = {3},
  number = {4},
  pages = {040305},
  issn = {2691-3399},
  doi = {10.1103/PRXQuantum.3.040305},
  urldate = {2026-04-24},
  langid = {english},
  author = {Dong, Y. and Lin, L. and Tong, Y.}
}

@misc{DP26a,
  title = {Constrained {{Optimal Polynomials}} for {{Quantum Linear System Solvers}} -- {{Numerical Data}}},
  year = 2026,
  month = apr,
  publisher = {Zenodo},
  doi = {10.5281/ZENODO.19694861},
  howpublished = {\url{https://zenodo.org/doi/10.5281/zenodo.19694860}},
  urldate = {2026-04-28},
  author = {Deiml, M. and Peterseim, D.}
}
